\documentclass{amsart}

\usepackage{graphicx}

\usepackage[textwidth=16cm, textheight=21cm]{geometry}
\usepackage{amsmath}
\usepackage{amssymb}
\usepackage{cite}

\newtheorem{theorem}{Theorem}

\newtheorem{proposition}[theorem]{Proposition}

\newtheorem{corollary}[theorem]{Corollary}
\theoremstyle{remark}
\newtheorem{example}[theorem]{Example}
\newtheorem{remark}[theorem]{Remark}

\newcommand{\arxiv}[1]{\href{http://www.arXiv.org/abs/#1}{arXiv:#1}}

\newcommand{\opl}{{\oplus}}

\newcommand{\rmin}{\mathbf{R}_{\min}}
\newcommand{\rmax}{\mathbf{R}_{\max}}
\newcommand{\suplim}{\sup\limits}
\newcommand{\sumlim}{\sum\limits}
\newcommand{\maxlim}{\max\limits}
\newcommand{\pd}[2]{\dfrac{\partial#1}{\partial#2}}
\newcommand{\CalB}{\mathcal{B}}

\newcommand{\maF}{\mathcal{F}}
\newcommand{\maD}{\mathcal{D}}

\newcommand{\cF}{{\mathcal F}}

\newcommand{\cN}{{\mathcal N}}

\newcommand{\0}{\mathbf{0}}
\newcommand{\1}{\mathbf{1}}
\newcommand{\cset}{\mathbf{C}}

\newcommand{\rset}{\mathbf{R}}
\newcommand{\maA}{\mathcal{A}}
\newcommand{\Log}{\mathop{\mathrm{Log}}}
\newcommand{\ovol}{\mathop{\mathrm{vol}}}

\def\C{\mathbf C}
\def\R{\mathbf R}
\def\maF{\mathcal F}
\def\maD{\mathcal D}

\def\cN{\mathcal N}
\def\Rmax{\rset_{\max}}

\newcommand{\rseth}{\widehat{\mathbf{R}}}
\newcommand{\rmaxh}{\rseth_{\max}}
\newcommand{\rminh}{\rseth_{\min}}
\newcommand{\smaxmin}{S_{\max,\min}}
\newcommand{\Mat}{\mathrm{Mat}}
\newcommand{\x}{\mathbf x}

\newcommand{\y}{\mathbf y}

\newcommand{\lx}{\underline{\x}}
\newcommand{\ly}{\underline{\y}}

\newcommand{\ux}{\overline{\x}}
\newcommand{\uy}{\overline{\y}}

\begin{document}

\fontsize{14pt}{20pt}
\selectfont

\title{Idempotent/tropical analysis, the Hamilton-Jacobi and Bellman equations}
\author{Grigory L. Litvinov}
\address{Grigory L. Litvinov, Institute for Information Transmission Problems,
B. Karetnyi per. 19/1, Moscow, 127994 Russia}
\email{glitvinov@gmail.com}
\thanks{To be published in Springer Lecture Notes in Mathematics.}

%
%

\begin{abstract}
Tropical and idempotent analysis with their relations to the
Hamilton-Jacobi and matrix Bellman equations are discussed. Some dequantization procedures are important in tropical and idempotent mathematics. In particular,
the Hamilton-Jacobi-Bellman equation is treated as a result of the Maslov
dequantization applied to the Schr\"{o}dinger equation. This leads to a linearity
of the Hamilton-Jacobi-Bellman equation over tropical algebras. The correspondence principle and the superposition principle of idempotent mathematics are formulated and examined. The matrix Bellman equation and its applications to optimization problems on graphs are
discussed. Universal algorithms for numerical algorithms in idempotent mathematics
are investigated. In particular, an idempotent version of interval analysis is briefly discussed.
\end{abstract}
\maketitle
\medskip

\begin{flushright}
{\em In dear memory of my beloved wife Irina.}
\end{flushright}

\section{Introduction}
\label{s:intro}
In these lecture notes we shall discuss some important problems of tropical and idempotent mathematics and especially those of idempotent and tropical analysis. Relations to the Hamilton-Jacobi and matrix Bellman equations will be examined. Applications of general principles of idempotent mathematics to numerical algorithms and their computer implementations will be discussed.

Tropical mathematics can be treated  as  a result of a
dequantization of the traditional mathematics as  the Planck
constant  tends to zero  taking imaginary values. This kind of
dequantization is known as the Maslov dequantization and it leads
to a mathematics over tropical algebras like the max-plus algebra.
The so-called idempotent dequantization is a generalization of the
Maslov dequantization. The idempotent dequantization leads to
mathematics over idempotent semirings (exact definitions see below
in sections 2 and 3). For  example,
the field of real or complex numbers can be treated as a quantum
object whereas idempotent semirings can be examined  as
"classical" or "semiclassical" objects  (a semiring is called
idempotent  if the  semiring addition is idempotent, i.e. $x
\oplus x = x$), see~\cite{Lit-07,LM-95,LM-96,LM-98}.
Some other dequantization procedures lead to interesting applications, e.g.,
to convex geometry, see below and \cite{LMS:07, LSz-05, LSz-07a}.

Tropical algebras are idempotent semirings (and semifields). Thus tropical mathematics is a part of idempotent mathematics. Tropical
algebraic geometry can be regarded as a result of the Maslov
dequantization applied to the traditional algebraic geometry (O.
Viro, G. Mikhalkin), see,
e.g.,~\cite{IMS:07,Mik-05,Mik-06,Vir-00,Vir-02,Vir-08}. There
are interesting relations and applications to the traditional
convex geometry.

    In the spirit of N.Bohr's correspondence principle there
is a (heuristic)  correspondence  between important, useful, and
interesting  constructions  and  results over fields and similar
constructions and results  over  idempotent  semirings.  A systematic application of
this correspondence  principle  leads  to  a variety  of theoretical
and applied results~\cite{Lit-07,LM-95,LM-96,LM-98,LM:05}, see Figure~\ref{f:bohr}.

The history of the subject is discussed, e.g., in~\cite{Lit-07}, with
extensive bibliography. See also
\cite{Car:79, CGQ-99, CGQ-04, CG:79, CG:95, LM-95, LM-96, LM-98,LMRS}.
\if{
\begin{figure}
\centering
\includegraphics[width=0.8\linewidth]{rand3by4}
\caption{Spectral function of \eqref{ABrand}}
\label{curioustest}
\end{figure}
}\fi

\begin{figure}[b]
\includegraphics[scale=1]{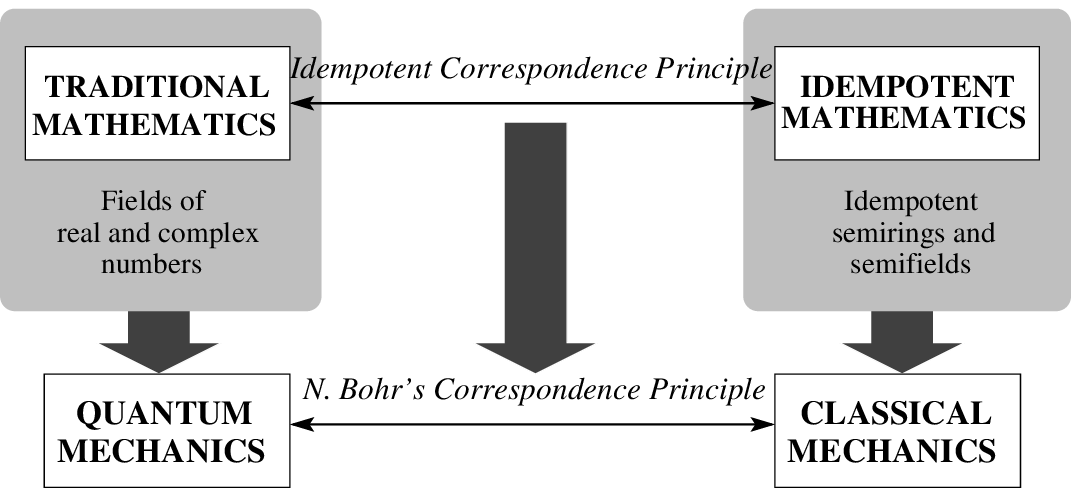}
%
%
\caption{Relations
between idempotent and traditional mathematics.}
\label{f:bohr}       
\end{figure}

V.P. Maslov's idempotent superposition principle means that many nonlinear problems
related to extremal problems are linear over suitable idempotent semirings. The principle is very important for applications including numerical and parallel
computations. See V.P. Maslov's original formulation in~\cite{Mas-86, Mas-87a, Mas-87b}, as well as~\cite{BCOQ, Car-71,
Car:79, CGQ-99, CGQ-04, CG:79, CG:95, KM:97, Lit-07,LM-95,LM-96,LM-98,LM:05,LMRS}, and below.

\section{The Maslov dequantization}
\label{s:dequant}
Let $\R$ and $\C$ be the fields of real and complex numbers. The
so-called max-plus algebra $\R_{\max}= \R\cup\{-\infty\}$ is
defined by the operations $x\oplus y=\max\{x, y\}$ and $x\odot y=
x+y$.

\begin{figure}[b]
\includegraphics[scale=0.5]{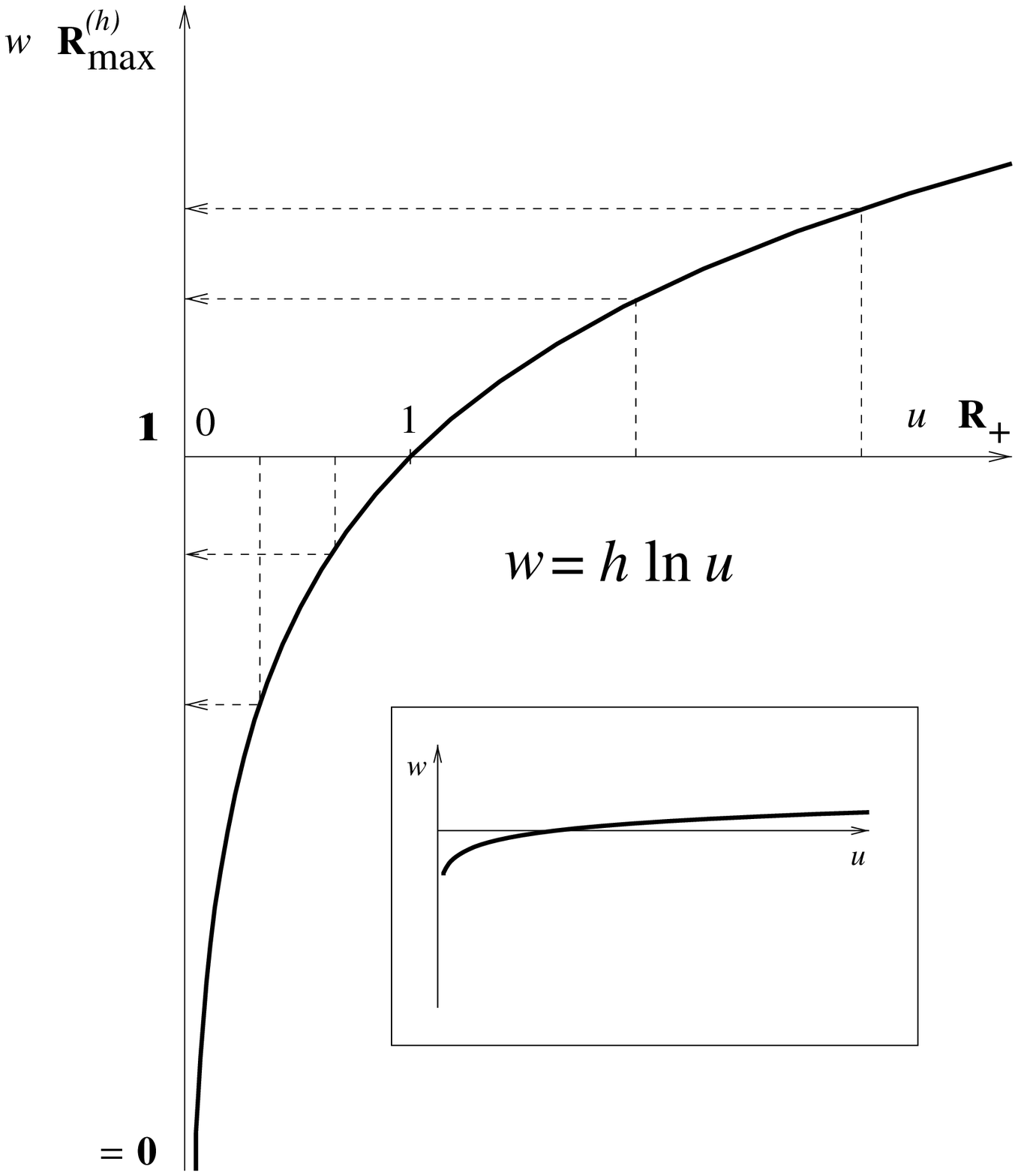}
%
%
\caption{Deformation of $\R_+$ to $\R^{(h)}$. Inset: the same for a
small value of $h$.}
\label{f:dequant}       
\end{figure}

The max-plus algebra can be seen as a result of the {\it Maslov
dequantization} of the semifield $\R_+$ of all nonnegative
numbers with the usual arithmetics. The change of variables
\begin{eqnarray*}
x\mapsto u=h\log x,
\end{eqnarray*}
where $h>0$, defines a map $\Phi_h\colon \R_+\to
\R\cup\{-\infty\}$, see Fig.~\ref{f:dequant}. Let the addition and multiplication
operations be mapped from \markboth{G.L. Litvinov}{Tropical Mathematics, Idempotent Analysis, Classical Mechanics, and Geometry} $\R_+$ to $\R\cup\{-\infty\}$ by $\Phi_h$, i.e.\
let
\begin{eqnarray*}
u\oplus_h v = h \log({\mbox{exp}}(u/h)+{\mbox{exp}}(v/h)),\quad u\odot v= u+ v,\\
\mathbf{0}=-\infty = \Phi_h(0),\quad \mathbf{1}= 0 = \Phi_h(1).
\end{eqnarray*}

\if{
\begin{figure}
\noindent\epsfig{file=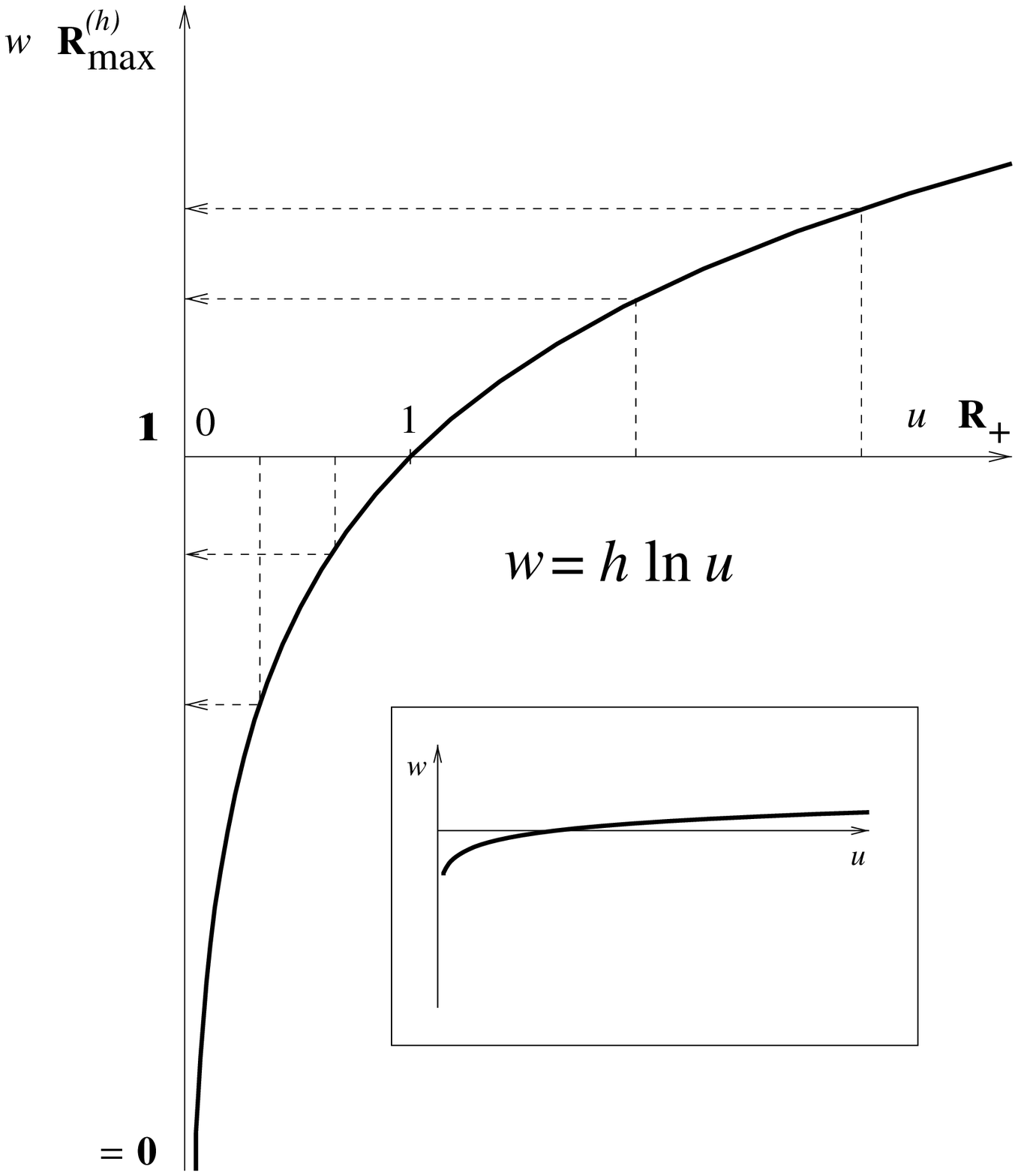,width=0.9\linewidth}
\vskip -3cm
\caption{Deformation of $\R_+$ to $\R^{(h)}$. Inset: the same for a
small value of $h$.}
\end{figure}
}\fi

It can be easily checked that $u\oplus_h v\to \max\{u, v\}$ as
$h\to 0$. This deformation of the algebraic
structure borrowed from~$\R_+$ brings us to the semifield $\R_{\max}$, known as the
{\em max-plus algebra}, with zero $\mathbf{0}= -\infty$ and unit
$\mathbf{1}=0$ .

The semifield $\R_{\max}$ is a typical example of an {\it
 idempotent semiring}; this is a semiring with idempotent addition, i.e.,
 $x\oplus x = x$ for arbitrary element
 $x$ of this semiring.

 The semifield $\R_{\max}$ is also called a \emph{tropical
 algebra}.The semifield $\R^{(h)}=\Phi_h(\R_+)$ with operations
 $\oplus_h$ and $\odot$ (i.e.$+$) is called a \emph{subtropical
 algebra}.

 The semifield $\R_{\min}=\R\cup\{+\infty\}$
with operations $\oplus={\min}$ and $\odot=+$
$(\mathbf{0}=+\infty, \mathbf{1}=0)$ is isomorphic to $\R_{\max}$.

The analogy with quantization is obvious; the parameter $h$ plays
the role of the Planck constant. The map $x\mapsto|x|$ and the
Maslov dequantization for $\R_+$ give us a natural transition from
the field $\C$ (or $\R$) to the max-plus algebra $\R_{\max}$. {\it
We will also call this transition the Maslov dequantization}. In fact
the Maslov dequantization corresponds to the usual Schr\"odinger
dequantization but for  imaginary values of the Planck constant (see below).
The transition from numerical fields to the max-plus algebra
$\R_{\max}$ (or similar semifields) in mathematical constructions
and results generates the so called {\it tropical mathematics}.
The so-called {\it idempotent dequantization} is a generalization
of the Maslov dequantization; this is the transition from basic fields to idempotent semirings in mathematical constructions and results without any deformation. The idempotent dequantization generates
the so-called \emph{idempotent mathematics}, i.e. mathematics over
idempotent semifields and semirings.

{\bf Remark.} The term 'tropical' appeared in~\cite{Sim-88} for a
discrete version of the max-plus algebra (as a suggestion of
Christian Choffrut). On the other hand V.P. Maslov used this term in
80s in his talks and works on economical applications of his
idempotent analysis (related to colonial politics). For the most
part of modern authors, 'tropical' means 'over $\R_{\max}$ (or
$\R_{\min}$)' and tropical algebras are $\R_{\max}$ and $\R_{\min}$.
The terms 'max-plus', 'max-algebra' and 'min-plus' are often used in
the same sense.

\section{Semirings and semifields. The idempotent correspondence principle}
\label{s:semi}

Consider a set $S$ equipped with two algebraic operations: {\it
addition} $\oplus$ and {\it multiplication} $\odot$. It is a {\it
semiring} if the following conditions are satisfied:
\begin{itemize}
\item the addition $\oplus$ and the multiplication $\odot$ are
associative; \item the addition $\oplus$ is commutative; \item the
multiplication $\odot$ is distributive with respect to the
addition $\oplus$:
\[x\odot(y\oplus z)=(x\odot y)\oplus(x\odot z)\]
and
\[(x\oplus y)\odot z=(x\odot z)\oplus(y\odot z)\]
for all $x,y,z\in S$.
\end{itemize}
A {\it unity} (we suppose that it exists) of a semiring $S$ is an element $\1\in S$ such that
$\1\odot x=x\odot\1=x$ for all $x\in S$. A {\it zero} (if it exists) of a
semiring $S$ is an element $\0\in S$ such that $\0\neq\1$ and
$\0\oplus x=x$, $\0\odot x=x\odot \0=\0$ for all $x\in S$. A
semiring $S$ is called an {\it idempotent semiring} if $x\oplus
x=x$ for all $x\in S$. A semiring $S$ with neutral element
$\1$ is called a {\it semifield} if every nonzero element of
$S$ is invertible with respect to the multiplication. For the theory
of semirings and semifields the reader is referred, e.g., to~\cite{Gol:99}.

The analogy with quantum physics discussed in Section~\ref{s:dequant} and
below leads to the following {\it idempotent correspondence principle}:

{\it There is a (heuristic) correspondence between important, useful and interesting
constructions and results over the field of complex (or real) numbers (or
the semifield of nonnegative numbers) and similar constructions and results over
idempotent semirings in the spirit of N.~Bohr's correspondence principle in
quantum theory~\cite{LM-95,LM-96,LM-98}.}

This principle can be also applied to algorithms and their software and hardware
implementations. Examples are discussed below; see also
\cite{Lit-07,LM-95,LM-96,LM-98, LMS-98, LMS-99,LMS-01,LMS-02,LS:09,LSz-02,
LSz-05,LSz-07a,LSz-07b}


\section{Idempotent analysis}
\label{s:idempan}
Idempotent analysis deals with functions taking their values in
an idempotent semiring and the corresponding function spaces.
Idempotent analysis was initially constructed by V.~P.~Maslov and
his collaborators and then developed by many authors. The subject
is presented in the book of V.~N.~Kolokoltsov and
V.~P.~Maslov~\cite{KM:97} (a version of this book in Russian was
published in 1994).

Let $S$ be an arbitrary semiring with idempotent addition $\oplus$
(which is always assumed to be commutative), multiplication
$\odot$, and unit $\1$. The set $S$ is equipped with
the {\it standard partial order\/}~$\preceq$: by definition, $a \preceq
b$ if and only if $a \oplus b = b$. If $S$ contains a zero element
$\0$, then all elements of $S$ are
nonnegative: $\0 \preceq$ $a$ for all $a \in S$. Due to the existence
of this order, idempotent analysis is closely related to the
lattice theory, theory of vector lattices, and theory of ordered
spaces. Moreover, this partial order allows to model a number of
basic ``topological'' concepts and results of idempotent analysis
on the purely algebraic level; this line of reasoning was examined
systematically in~\cite{Lit-07}--~\cite{LSz-07b}
and~\cite{CGQ-04}.

Calculus deals mainly with functions whose values are numbers. The
idempotent analog of a numerical function is a map $X \to S$,
where $X$ is an arbitrary set and $S$ is an idempotent semiring.
Functions with values in $S$ can be added, multiplied by each
other, and multiplied by elements of $S$ pointwise.

The idempotent analog of a linear functional space is a set of
$S$-valued functions that is closed under addition of functions
and multiplication of functions by elements of $S$, or an
$S$-semimodule. Consider, e.g., the $S$-semimodule $B(X, S)$ of
all functions $X \to S$ that are bounded in the sense of the
standard order on $S$.

If $S = \rmax$, then the idempotent analog of integration is
defined by the formula
\begin{equation}
\label{id-integral}
I(\varphi) = \int_X^{\oplus} \varphi (x)\, dx     = \sup_{x\in X}
\varphi (x),
\end{equation}
where $\varphi \in B(X, S)$. Indeed, a Riemann sum of the form
$\sumlim_i \varphi(x_i) \cdot \sigma_i$ corresponds to the
expression $\bigoplus\limits_i \varphi(x_i) \odot \sigma_i =
\maxlim_i \{\varphi(x_i) + \sigma_i\}$, which tends to the
right-hand side of~\eqref{id-integral} as $\sigma_i \to 0$. Of course, this is a
purely heuristic argument.

Formula~\eqref{id-integral} defines the \emph{idempotent} (or \emph{Maslov})
\emph{integral} not only for functions taking values in $\rmax$,
but also in the general case when any of bounded (from above)
subsets of~$S$ has the least upper bound.

An \emph{idempotent} (or \emph{Maslov}) \emph{measure} on $X$ is
defined by the formula $m_{\psi}(Y) = \suplim_{x \in Y} \psi(x)$, where $\psi
\in B(X,S)$ is a fixed function. The integral with respect to this measure is defined
by the formula
\begin{equation}
\label{id-measure}
I_{\psi}(\varphi)
    = \int^{\oplus}_X \varphi(x)\, dm_{\psi}
    = \int_X^{\oplus} \varphi(x) \odot \psi(x)\, dx
    = \sup_{x\in X} (\varphi (x) \odot \psi(x)).
\end{equation}

Obviously, if $S = \rmin$, then the standard order is
opposite to the conventional order $\leq$, so in this case
equation~\eqref{id-measure} takes the form
\begin{equation*}
   \int^{\oplus}_X \varphi(x)\, dm_{\psi}
    = \int_X^{\oplus} \varphi(x) \odot \psi(x)\, dx
    = \inf_{x\in X} (\varphi (x) \odot \psi(x)),
\end{equation*}
where $\inf$ is understood in the sense of the conventional order
$\leq$.

We shall see that in idempotent analysis measures and generalized functions (versions of distributions in the sense of L.~Schwartz) are generated by usual
functions. For example the $\delta$-functional $\delta_y:\ \varphi(\cdot)\mapsto\varphi(y)$
is generated by the function
\begin{equation*}
\delta_y(x)=
\begin{cases}
\1, &\text{if $x=y$},\\
\0, &\text{if $x\neq y$}.
\end{cases}
\end{equation*}
It is clear that $$\varphi(y)=\int^{\oplus}_X \delta_y(x)\odot\varphi(x) dx=
\sup_x (\delta_y(x)\odot\varphi(x))$$ .

\section{The superposition principle and linear equations}
\label{s:super}
\subsection{Heuristics}

Basic equations of quantum theory are linear; this is the
superposition principle in quantum mechanics. The Hamilton--Jacobi
equation, the basic equation of classical mechanics, is nonlinear
in the conventional sense. However, it is linear over the
semirings $\rmax$ and $\rmin$. Similarly, different versions of
the Bellman equation, the basic equation of optimization theory,
are linear over suitable idempotent semirings; this is
V.~P.~Maslov's idempotent superposition principle, see
\cite{Mas-86,Mas-87a,Mas-87b}. More generally, the idempotent superposition principle means
that although some important problems and equations (related to extremal problems, e.g., optimization problems, the Bellman equation and its instances, the
Hamilton-Jacobi equation) are nonlinear in the usual sense, they can be treated as linear over appropriate idempotent semirings.
For instance, the finite-dimensional
stationary Bellman equation can be written in the form $X = H
\odot X \oplus F$, where $X$, $H$, $F$ are matrices with
coefficients in an idempotent semiring $S$ and the unknown matrix
$X$ is determined by $H$ and $F$, see below and~\cite{Car-71, Car:79, BCOQ, CG:79, CG:95, GM:79, GM:10}. In
particular, standard problems of dynamic programming and the
well-known shortest path problem correspond to the cases $S =
\rmax$ and $S =\rmin$, respectively. It is known that principal
optimization algorithms for finite graphs correspond to standard
methods for solving systems of linear equations of this type
(i.e., over semirings). Specifically, Bellman's shortest path
algorithm corresponds to a version of Jacobi's algorithm, Ford's
algorithm corresponds to the Gauss--Seidel iterative scheme,
etc.~\cite{Car-71, Car:79}.

The linearity of the Hamilton--Jacobi equation over $\rmin$ and
$\rmax$, which is the result of the Maslov dequantization of the
Schr{\"o}\-din\-ger equation, is closely related to the
(conventional) linearity of the Schr{\"o}\-din\-ger equation and
can be deduced from this linearity. Thus, it is possible to borrow
standard ideas and methods of linear analysis and apply them to a
new area.

Consider a classical dynamical system specified by the Hamiltonian
$$
   H = H(p,x) = \sum_{i=1}^N \frac{p^2_i}{2m_i} + V(x),
$$
where $x = (x_1, \dots, x_N)$ are generalized coordinates, $p =
(p_1, \dots, p_N)$ are generalized momenta, $m_i$ are generalized
masses, and $V(x)$ is the potential. In this case the Lagrangian
$L(x, \dot x, t)$ has the form
$$
   L(x, \dot x, t)
    = \sum^N_{i=1} m_i \frac{\dot x_i^2}2 - V(x),
$$
where $\dot x = (\dot x_1, \dots, \dot x_N)$, $\dot x_i = dx_i /
dt$. The value function $S(x,t)$ of the action functional has the
form
$$
   S = \int^t_{t_0} L(x(t), \dot x(t), t)\, dt,
$$
where the integration is performed along the actual trajectory of
the system.  The classical equations of motion are derived as the
stationarity conditions for the action functional (the Hamilton
principle, or the least action principle).

For fixed values of $t$ and $t_0$ and arbitrary trajectories
$x(t)$, the action functional $S=S(x(t))$ can be considered as a
function taking the set of curves (trajectories) to the set of
real numbers which can be treated as elements of  $\rmin$. In this
case the minimum of the action functional can be viewed as the
Maslov integral of this function over the set of trajectories or
an idempotent analog of the Euclidean version of the Feynman path
integral. The minimum of the action functional corresponds to the
maximum of $e^{-S}$, i.e. idempotent integral
$\int^{\oplus}_{\{paths\}} e^{-S(x(t))} D\{x(t)\}$ with respect to
the max-plus algebra $\rset_{\max}$. Thus the least action
principle can be considered as an idempotent version of the
well-known Feynman approach to quantum mechanics.  The
representation of a solution to the Schr{\"o}\-din\-ger equation
in terms of the Feynman integral corresponds to the
Lax--Ole\u{\i}nik solution formula for the Hamilton--Jacobi
equation.

Since $\partial S/\partial x_i = p_i$, $\partial S/\partial t =
-H(p,x)$, the following Hamilton--Jacobi equation holds:
\begin{equation}
\label{e:HJ}
\pd{S}{t} + H \left(\pd{S}{x_i}, x_i\right)= 0.
\end{equation}

Quantization leads to the Schr\"odinger equation
\begin{equation}
\label{e:Sch}
  -\frac{\hbar}i \pd{\psi}{t}= \widehat H \psi = H(\hat p_i, \hat x_i)\psi,
\end{equation}
where $\psi = \psi(x,t)$ is the wave function, i.e., a
time-dependent element of the Hilbert space $L^2(\rset^N)$, and
$\widehat H$ is the energy operator obtained by substitution of
the momentum operators $\widehat p_i = {\hbar \over i}{\partial
\over \partial x_i}$ and the coordinate operators $\widehat x_i
\colon \psi \mapsto x_i\psi$ for the variables $p_i$ and $x_i$ in
the Hamiltonian function, respectively. This equation is linear in
the conventional sense (the quantum superposition principle). The
standard procedure of limit transition from the Schr\"odinger
equation to the Hamilton--Jacobi equation is to use the following
ansatz for the wave function:  $\psi(x,t) = a(x,t)
e^{iS(x,t)/\hbar}$, and to keep only the leading order as $\hbar
\to 0$ (the `semiclassical' limit).

Instead of doing this, we switch to imaginary values of the Planck
constant $\hbar$ by the substitution $h = i\hbar$, assuming $h >
0$. Then the Schr\"odinger equation~\eqref{e:Sch} becomes similar to the
heat equation:
\begin{equation}
\label{e:heat}
   h\pd{u}{t} = H\left(-h\frac{\partial}{\partial x_i}, \hat x_i\right) u,
\end{equation}
where the real-valued function $u$ corresponds to the wave
function $\psi$. A similar idea (a switch to imaginary time) is
used in the Euclidean quantum field theory; let us remember that
time and energy are dual quantities.

Linearity of equation~\eqref{e:Sch} implies linearity of equation~\eqref{e:heat}. Thus if $u_1$ and $u_2$ are solutions of~\eqref{e:heat}, then so is their linear
combination
\begin{equation}
\label{heat-linear}
   u = \lambda_1 u_1 + \lambda_2 u_2.
\end{equation}

Let $S = h \ln u$ or $u = e^{S/h}$ as in Section~\ref{s:dequant} above. It can
easily be checked that equation~\eqref{e:heat} thus turns to
\begin{equation}
\label{e:burgers}
 \pd{S}{t}= V(x) + \sum^N_{i=1} \frac1{2m_i}\left(\pd{S}{x_i}\right)^2
    + h\sum^n_{i=1}\frac{1}{2m_i}\frac{\partial^2 S}{\partial x^2_i}.
\end{equation}
Thus we have a transition from~\eqref{e:HJ} to~\eqref{e:burgers} by means of the change of
variables $\psi = e^{S/h}$. Note that $|\psi| = e^{ReS/h}$ , where
Re$S$ is the real part  of $S$. Now let us consider $S$ as a real
variable. The equation~\eqref{e:burgers} is nonlinear in the conventional sense.
However, if $S_1$ and $S_2$ are its solutions, then so is the
function
\begin{equation}
\label{burgers-linear}
S = \lambda_1 \odot S_1 \opl_h \lambda_2\odot S_2
\end{equation}
obtained from~\eqref{heat-linear} by means of the substitution
$S = h \ln u$. Here
the generalized multiplication $\odot$ coincides with the ordinary
addition and the generalized addition $\opl_h$ is the image of the
conventional addition under the above change of variables.  As $h
\to 0$, we obtain the operations of the idempotent semiring
$\rmax$, i.e., $\oplus = \max$ and $\odot = +$, and equation~\eqref{e:burgers}
becomes the Hamilton--Jacobi equation~\eqref{e:HJ}, since the third term
in the right-hand side of equation~\eqref{e:burgers} vanishes.

Thus it is natural to consider the limit function $S = \lambda_1
\odot S_1 \oplus \lambda_2 \odot S_2$ as a solution of the
Hamilton--Jacobi equation and to expect that this equation can be
treated as linear over $\rmax$. This argument (clearly, a
heuristic one) can be extended to equations of a more general
form. For a rigorous treatment of (semiring) linearity for these
equations see, e.g., \cite{KM:97,LM:05,Roub}. Notice that if
$h$ is changed to $-h$, then we have that the resulting
Hamilton--Jacobi equation is linear over $\rmin$.

The idempotent superposition principle indicates that there exist
important nonlinear (in the traditional sense) problems that are
linear over idempotent semirings. The idempotent linear functional
analysis (see below) is a natural tool for investigation of those
nonlinear infinite-dimensional problems that possess this
property.

\subsection{The Cauchy problem for the Hamilton-Jacobi equations}

A rigorous ``idempotent'' appproach to the investigation of the Hamilton-Jacobi
equation was developed by V.N. Kolokoltsov and V.P. Maslov~\cite{KM:97}
(a Russian version of this book was published in 1994); see
also~\cite{McE:10, Roub, Sub:95, Sub-96}.

Let us consider, inspired by a long tradition, the well-known Cauchy problem
for the Hamilton-Jacobi equation~\eqref{e:HJ}. Given the action function at
time $T$
\begin{equation}
\label{cauchy}
S(T,x)=S_T(x)=\varphi(x),\quad x\in\R^N,
\end{equation}
the Cauchy problem asks to reconstruct $S(t,x)$ for $x\in\R^N$ during the time interval $0\leq t\leq T$.

We shall discuss the min-plus linearity of this problem and denote by $U_t$ the resolving operator, i.e. the map which assigns to each given $S_T(x)$ the
solution $S(t,x)$ of the Cauchy problem in the interval $0\leq t\leq T$. Then the
map $U_t$, for each $t$, is a linear (over $\rmin$) operator in the space
LSC$(\R^n,\rmin)$ of lower semicontinuous functions taking their values in $\rmin$.
Moreover $U_t$ is an integral operator (in the sense of idempotent mathematics) of the form:
\begin{equation}
(U_t\varphi)(x)=\int^{\oplus}\varphi(y) K_t(x,y) dy=\inf_y \{\varphi(y)+K_t(x,y)\},
\end{equation}
where $K_t(x,y)$, as a function of $y\in\R^n$, is bounded from below and lower semicontinuous. See~\cite{KM:97,Roub} for details.

The operator $U_t$ (as well as other integral operators, see Section~\ref{s:funcan}
below) has the following property:
\begin{equation}
U_t(\bigoplus_{\nu} \varphi_{\nu})=\bigoplus_{\nu} (U_t \varphi_{\nu}),
\end{equation}
where $\{\varphi_{\nu}\}$ is a bounded set of elements in LSC$(\R^n,\rmin)$. So if
we have such a family of functions $S_{\nu}(T,x)$ and
$S(T,x)=\int^{\oplus} S_{\nu}(T,x) d\nu=\inf_{\nu}(S_{\nu}(T,x))$, then the solution
of the Cauchy problem is expressed as $S(t,x)=\inf_{\nu}(S_{\nu}(t,x))$.

Relations between the ``idempotent approach'', viscosity solutions and minimax
solutions in the sense of A.I.~Subbotin~\cite{Sub:95,Sub-96} are examined, e.g.,
in~\cite{Roub} in details; see also W.M.~McEneaney~\cite{McE:10}. To this end,
let us mention that more general
Hamiltonians of the form $H=H(t,x,p)$ (satisfying some additional conditions) and
different kinds of solution spaces are also considered in the literature.

The situation is similar for the Cauchy problem for the homogeneous Hamilton-Jacobi
equation
\begin{equation*}
\pd{S}{t}+H(\pd{S}{x})=0,\quad S_{t=0}=S_0(x),
\end{equation*}
where $H:\R^n\mapsto \R$ is a convex (not strictly) first order homogeneous
function
\begin{equation*}
H(p)=\sup\limits_{(f,g)\in V} (f\cdot p+g),\quad f\in\R^n,\quad g\in\R,
\end{equation*}
and $V$ is a compact set in $\R^{n+1}$. See~\cite{KM:97}.

To develop a rigorous ``idempotent'' approach to differential equations and other
problems, one needs an idempotent version of analysis and, especially, functional analysis.
See Section~\ref{s:funcan} below.

\section{Convolution and the Fourier--Legendre transform}
\label{s:fourier-legendre}

Let $G$ be a group. Then the space $\CalB(G, \rset_{\max})$ of all
bounded functions $G\to\rset_{\max}$ (see above) is an idempotent
semiring with respect to the following analog $\circledast$ of the
usual convolution:
$$
   (\varphi(x)\circledast\psi)(g)=
    = \int_G^{\oplus} \varphi (x)\odot\psi(x^{-1}\cdot g)\, dx=
\sup_{x\in G}(\varphi(x)+\psi(x^{-1}\cdot g)).
$$
Of course, it is possible to consider other ``function spaces''
(and other basic semirings instead of $\rset_{\max}$).

Let $G=\rset^n$, where $\rset^n$ is considered as a topological
group with respect to the vector addition. The conventional
Fourier--Laplace transform is defined as
\begin{equation}
\label{e:FL}
   \varphi(x) \mapsto \tilde{\varphi}(\xi)
    = \int_G e^{i\xi \cdot x} \varphi (x)\, dx
\end{equation}
where $e^{i\xi \cdot x}$ is a character of the group $G$, i.e., a
solution of the following functional equation:
$$
   f(x + y) = f(x)f(y).
$$
The idempotent analog of this equation is
$$
   f(x + y) = f(x) \odot f(y) = f(x) + f(y),
$$
so ``continuous idempotent characters'' are linear functionals of
the form $x \mapsto \xi \cdot x = \xi_1 x_1 + \dots + \xi_n x_n$.
As a result, the transform in~\eqref{e:FL} assumes the form
\begin{equation}
\label{e:Legendre}
   \varphi(x) \mapsto \tilde{\varphi}(\xi)
    = \int_G^\oplus \xi \cdot x \odot \varphi (x)\, dx
   = \sup_{x\in G} (\xi \cdot x + \varphi (x)).
\end{equation}
The transform in~\eqref{e:Legendre} is the {\it Legendre transform\/}
(up to some change of notation) \cite{Mas-87b}; transforms of this kind
establish the correspondence between the Lagrangian and the
Hamiltonian formulations of classical mechanics. The Legendre
transform generates an idempotent version of harmonic analysis for
the space of convex functions, see, e.g., \cite{MT:03}.

Of course, this construction can be generalized to different
classes of groups and semirings. Transformations of this type
convert the generalized convolution $\circledast$ to the pointwise
(generalized) multiplication and possess analogs of some important
properties of the usual Fourier transform.

The examples discussed in this sections can be treated as
fragments of an idempotent version of the representation theory,
see, e.g., \cite{LMS-02}. In particular, ``idempotent''
representations of groups can be examined as representations of
the corresponding convolution semirings (i.e. idempotent group
semirings) in semimodules.

\section{Idempotent functional analysis}
\label{s:funcan}

Many other idempotent analogs may be given, in particular, for
basic constructions and theorems of functional analysis.
Idempotent functional analysis is an abstract version of
idempotent analysis. For the sake of simplicity take $S=\rmax$ and
let $X$ be an arbitrary set. The idempotent integration can be
defined by the formula (1), see above. The functional $I(\varphi)$
is linear over $S$ and its values correspond to limiting values of
the corresponding analogs of Lebesgue (or Riemann) sums. An
idempotent scalar product of functions $\varphi$ and $\psi$ is
defined by the formula
$$
\langle\varphi,\psi\rangle = \int^{\oplus}_X
\varphi(x)\odot\psi(x)\, dx = \sup_{x\in
X}(\varphi(x)\odot\psi(x)).
$$
So it is natural to construct idempotent analogs of integral
operators in the form
\begin{equation}
\label{e:int-operator}
\varphi(y) \mapsto (K\varphi)(x) = \int^{\oplus}_Y K(x,y)\odot\varphi(y)\, dy = \sup_{y\in Y}\{K(x,y)+\varphi(y)\},
\end{equation}
where $\varphi(y)$ is an element of a space of functions defined
on a set $Y$, and $K(x,y)$ is an $S$-valued function on $X\times
Y$. Of course, expressions of this type are standard in
optimization problems.

Recall that the definitions and constructions described above can
be extended to the case of idempotent semirings which are
conditionally complete in the sense of the standard order. Using
the Maslov integration, one can construct various function spaces
as well as idempotent versions of the theory of generalized
functions (distributions). For some concrete idempotent function
spaces it was proved that every `good' linear operator (in the
idempotent sense) can be presented in the form~\eqref{e:int-operator}; this is an
idempotent version of the kernel theorem of L.~Schwartz; results
of this type were proved by V.~N.~Kolokoltsov, P.~S.~Dudnikov and
S.~N.~Samborski\u\i, I.~Singer, M.~A.~Shubin and others. So every
`good' linear functional can be presented in the form
$\varphi\mapsto\langle\varphi,\psi\rangle$, where
$\langle,\rangle$ is an idempotent scalar product.

In the framework of idempotent functional analysis results of this
type can be proved in a very general situation. In \cite{LMS-98,
LMS-99,LMS-01,LMS-02,LSz-02,LSz-07b} an algebraic
version of the idempotent functional analysis is developed; this
means that basic (topological) notions and results are simulated
in purely algebraic terms (see below). The treatment covers the subject  from
basic concepts and results (e.g., idempotent analogs of the
well-known theorems of Hahn-Banach, Riesz, and Riesz-Fisher) to
idempotent analogs of A.~Grothendieck's concepts and results on
topological tensor products, nuclear spaces and operators.
Abstract idempotent versions of the kernel theorem are formulated. Note that
the transition from the usual theory to idempotent functional
analysis may be very nontrivial; for example, there are many
non-isomorphic idempotent Hilbert spaces. Important results on
idempotent functional analysis (duality and separation theorems)
were obtained by G.~Cohen, S.~Gaubert, and J.-P.~Quadrat.
Idempotent functional analysis has received much attention in the
last years, see,
e.g.,~\cite{AGK, CGQ-04, GM:79, GM:10, Gun:98, MS:92, Shub-92},~\cite{KM:97}--~\cite{LSz-07b} and works
cited in~\cite{Lit-07}. All the results presented in this section are proved in \cite{LMS-01}  (Subsections~\ref{ss:semimod} -- \ref{ss:linoper}) and in \cite{LSz-07b} (Subsections~\ref{ss:funcsem} -- \ref{ss:intrepr})

\subsection{Idempotent semimodules and idempotent linear spaces}

\label{ss:semimod}

An additive semigroup $S$ with commutative addition
$\oplus$ is called an {\it idempotent semigroup} if the relation
$x\oplus x=x$ is fulfilled for all elements $x\in S$. If $S$
contains a neutral element, this element is denoted by the
symbol $\0$. Any idempotent semigroup is a partially ordered set with respect
to the following standard order: $x\preceq y$ if and only
if $x\oplus y=y$. It is obvious that this order is well defined
and $x\oplus y=\sup \{x, y\}$. Thus, any idempotent semigroup is an upper semilattice;
moreover, the concepts of idempotent semigroup and upper semilattice coincide, see \cite{Bir:67}. An idempotent semigroup $S$ is called $a$-{\it complete}
(or {\it algebraically complete}) if it is complete as an ordered set, i.e., if any subset $X$ in $S$  has the least upper bound
$\sup(X)$ denoted by $\oplus X$ and the greatest lower bound $\inf (X)$ denoted by $\wedge X$. This semigroup is called $b$-{\it complete}
(or {\it boundedly complete}), if any bounded above subset $X$ of this
semigroup (including the empty subset) has the least upper bound
$\oplus X$ (in this case, any nonempty subset $Y$ in $S$
has the greatest lower bound $\wedge Y$ and $S$ in a lattice). Note that
any $a$-complete or $b$-complete idempotent semiring has the zero element $\0$
that coincides with 
$\oplus{\emptyset}$,
where 
$\emptyset$ is the empty set. Certainly,
$a$-completeness implies the $b$-completeness.
Completion by means of cuts \cite{Bir:67} yields an embedding
$S\to \widehat S$ of an arbitrary idempotent semigroup $S$ into an $a$-complete idempotent semigroup
$\widehat S$ (which is called a {\it normal completion} of $S$); in addition,
$\widehat{\widehat S} = S$. The $b$-completion procedure
$S \to \widehat S_b$ is defined similarly: if $S \ni \infty =\sup S$,
then $\widehat S_b$ =$\widehat S$; otherwise,
$\widehat S =\widehat S_b \cup \{ \infty \}$.
An arbitrary $b$-complete idempotent semigroup  $S$ also may differ from $\widehat S$ only
by the element $\infty =\sup S$.

Let $S$ and $T$ be $b$-complete idempotent semigroups.
Then, a homomorphism $f : S\to T$ is said to be a $b$-{\it homomorphism} if
 $f (\oplus X) = \oplus f(X)$ for any bounded subset $X$ in $S$.
If the  $b$-homomorphism $f$ is extended to a homomorphism
$\widehat S\to \widehat T$ of the correesponding normal completions and
$f(\oplus X) = \oplus f(X)$ for all $X\subset S$, then $f$
is said to be an $a$-{\it homomrphism}. An idempotent semigroup $S$ equipped with a topology such
that the set $\{ s\in S\vert s\preceq b\}$ is closed in this topology
for any $b\in S$ is called a {\it topological idempotent semigroup} $S$.

\begin {proposition} Let $S$ be an $a$-complete topological
idempotent semigroup and $T$ be a  $b$-complete topological idempotent semigroup such that, for
any nonempty subsemigroup $X$ in $T$, the element $\oplus X$ is
contained in the topological closure of $X$ in $T$. Then, a
homomorphism $f : T\to S$ that maps zero into zero
is an $a$-homomorphism if and only
if the mapping $f$ is lower semicontinuous in the sense that
the set $\{ t\in T\vert f(t)\preceq s\}$ is closed in $T$ for any $s\in S$.
\end{proposition}
\medskip
An idempotent semiring $K$ is called $a$-{\it complete} (respectively
$b$-{\it complete})
if $K$ is an $a$-complete (respectively $b$-complete) idempotent semigroup and,
for any subset (respectively, for any bounded subset) $X$ in $K$
and any $k\in K$, the generalized distributive laws
$k\odot(\oplus X)=\oplus (k\odot X)$ and $(\oplus X)\odot k =
\oplus(X\odot k)$ are fulfilled. Generalized distributivity
implies that any $a$-complete or $b$-complete idempotent semiring has a zero
element that coincides with 
$\oplus{\emptyset}$,
where 
$\emptyset$
is the empty set.

The set $\R(\max, +)$ of real numbers equipped with the idempotent
addition $\oplus=\max$ and multiplication $\odot=+$ is  an
idempotent semiring; in this case, $\1 = 0$. Adding the element $\0=-\infty$
to this semiring, we obtain a $b$-complete semiring
$\R_{\max} = \R\cup \{-\infty\}$
with the same operations and the zero element.
Adding the element $+\infty$ to $\R_{\max}$ and asumming that
$\0\odot (+\infty)=\0$ and $x\odot (+\infty) = +\infty$
for $x\neq\0$ and $x\oplus (+\infty) = +\infty$ for any $x$, we
obtain the $a$-complete idempotent semiring $\widehat{\R}_{\max} = \R_{\max}\cup \{+\infty\}$.
The standard order on $\R(\max, +)$, $\R_{\max}$ and
$\widehat{\R}_{\max}$ coincides with the ordinary order.
The semirings $\R(\max,+)$ and  $\R_{\max}$ are semifields.
On the contrary, an $a$-complete semiring that does not coincide with
$\{\0, \1\}$ cannot be a semifield. An important class
of examples is related to
(topological) vector lattices (see, for example, \cite{Bir:67} and
\cite{Scha}, Chapter 5). Defining the sum $x\oplus y$ as $\sup\{x, y\}$ and
the multiplication $\odot$ as the addition of vectors, we can
interpret the vector lattices as idempotent semifields.
Adding the zero element $\0$ to a complete vector lattice
(in the sense of \cite{Bir:67, Scha}), we obtain a $b$-complete semifield.
If, in addition, we add the infinite element, we obtain an
$a$-complete idempotent semiring (which, as an ordered set, coincides with the normal
completion of the original lattice).

{\bf Important definitions.} Let $V$ be an idempotent semigroup and  $K$ be an
idempotent semiring. Suppose that a multiplication $k, x\mapsto k\odot x$
of all elements from $K$ by the elements from
$V$ is defined; moreover, this multiplication is associative
and distributive with respect
to the addition in $V$ and $\1\odot x = x$, $\0\odot x=\0$
for all $x\in V$. In this case, the semigroup $V$ is called an
{\it idempotent semimodule} (or simply, a {\it semimodule})
 over $K$. The element
 $\0_V\in V$ is called the {\it zero} of the semimodule $V$ if $k\odot\0_V=\0_V$
and $\0_V\oplus x = x$ for any $k\in K$ and $x\in V$.
Let $V$ be a semimodule over a $b$-complete idempotent semiring $K$.
This semimodule is called $b$-{\it complete} if it is $b$-complete as
an idempotent semiring and, for any bounded subsets $Q$ in $K$ and $X$ in  $V$,
the generalized distributive laws $(\oplus Q)\odot x = \oplus (Q\odot x)$ and
$k\odot (\oplus X) = \oplus (k\odot X)$ are fulfilled for all
$k\in K$ and $x\in X$. This semimodule is called $a$-{\it complete}
if it is $b$-complete and contains the element $\infty = \sup V$.

A semimodule $V$ over a $b$-complete semifield $K$ is said to be
an {\it idempotent} $a$-{\it space} ($b$-{\it space})
 if this semimodule is
$a$-complete (respectively, $b$-complete) and the equality
$(\wedge Q)\odot x = \wedge (Q\odot x)$ holds for any nonempty subset
$Q$ in $K$ and any $x\in V$, $x\neq \infty = \sup V$.
The normal completion $\widehat V$ of a $b$-space $V$ (as an idempotent semigroup) has the
structure of an idempotent $a$-space (and may differ from $V$ only
by the element $\infty= \sup V$).

Let $V$ and $W$ be idempotent semimodules over an idempotent semiring
$K$.  A mapping $p: V\to W$ is said to be {\it linear} (over $K$) if
$$
p(x\oplus y)=p(x)\oplus p(y) \mbox{ and } p(k\odot x)=k\odot p(x)
$$
for any $x, y\in V$ and $k\in K$. Let the semimodules $V$ and $W$ be
 $b$-complete. A linear mapping $p: V\to W$ is said to be $b$-{\it linear}
if it is a $b$-homomorphism of the idempotent semigroup; this mapping is said to be
$a$-{\it linear} if it can be extended to an $a$-homomorphism of the normal
completions $\widehat V$ and $\widehat W$. Proposition 7.1 (see above)
shows that $a$-linearity simulates (semi)continuity for
linear mappings. The normal completion $\widehat K$ of the
semifield $K$ is a semimodule over $K$. If $W= \widehat K$, then the
linear mapping $p$ is called a {\it linear functional}.

Linear, $a$-linear and $b$-linear mappings are also called {\it linear, a-linear} and {\it b-linear operators} respectively.

Examples of idempotent semimodules and spaces that are the most important
for analysis are either subsemimodules of topological vector lattices
\cite{Scha} (or coincide with them) or are dual to them, i.e., consist
of linear functionals subject to some regularity condition,
for example, consist of $a$-linear functionals. Concrete examples of
idempotent semimodules and spaces of functions (including spaces of
bounded, continuous, semicontinuous, convex, concave  and Lipschitz functions)
see in \cite{KM:97, LMS-99, LMS-01, LSz-07b} and below.

\subsection{Basic results}
\label{ss:basic}
Let $V$ be an idempotent $b$-space over  a
$b$-complete semifield $K$,  $x\in \widehat V$. Denote by
$x^*$ the functional $V\to \widehat K$ defined by the formula
$x^* (y)=\wedge \{ k\in K \vert y\preceq k\odot x\}$,
where $y$ is an arbitrary fixed element from $V$.

\begin{theorem}
\label{t:xfunc}
For any $x\in \widehat V$ the functional $x^*$
is $a$-linear. Any nonzero $a$-linear functional $f$ on $V$
is given by $f = x^*$ for a unique suitable element $x\in V$.
If $K\neq \{\0, \1\}$, then $x=\oplus\{ y\in V\vert f(y)\preceq\1\}$.
\end{theorem}

Note that results of this type obtained earlier concerning the structure of
linear functionals cannot be carried over to subspaces and
subsemimodules.

A subsemigroup $W$ in $V$ closed with respect to the
multiplication by an arbitrary element from $K$ is called a
$b$-{\it subspace} in $V$ if the imbedding $W\to V$ can be extended to a
$b$-linear mapping. The following result is obtained from Theorem~\ref{t:xfunc}
and is the idempotent version of the Hahn--Banach theorem.

\begin{theorem} Any $a$-linear functional defined on a
$b$-subspace $W$ in $V$ can be extended to an $a$-linear functional
on $V$. If $x, y\in V$ and $x\neq y$, then there exists an
$a$-linear functional $f$ on $V$ that separates the elements
$x$ and $y$, i.e., $f(x)\neq f(y)$.
\end{theorem}

The following statements are easily derived from the definitions and
can be regarded as the analogs of the well-known results of the traditional functional analysis (the Banach--Steinhaus and the closed-graph theorems).

\begin{proposition}
\label{p:BS1}
 Suppose that $P$ is a family of $a$-linear
mappings of an $a$-space $V$ into an $a$-space $W$ and the mapping
$p : V\to W$ is the pointwise sum of the mappings of this family,
i.e., $p(x) =\sup \{ p_{\alpha} (x)\vert p_{\alpha} \in P\}$.
Then the mapping $p$ is $a$-linear.
\end{proposition}

\begin{proposition} Let $V$ and $W$ be  $a$-spaces. A linear
mapping $p : V \to W$ is $a$-linear if and only if its graph
$\Gamma$ in $V\times W$ is closed with respect to passing
to sums
(i.e., to  least upper bounds)
of its arbitrary subsets.
\end{proposition}

In \cite{CGQ-04} the basic results were generalized for the
case of semimodules over the so-called reflexive $b$-complete semirings.

\subsection{Idempotent $b$-semialgebras}

Let $K$ be a $b$-complete semifield and $A$ be an
idempotent $b$-space over $K$ equipped with the structure of a semiring
compatible with the multiplication $K\times A\to A$ so that
the associativity of the multiplication is preserved. In this case,
$A$ is called an {\it idempotent $b$-semialgebra} over~$K$.

\begin{proposition}
\label{p:semialg}
For any invertible element $x\in A$ from the
$b$-semialgebra $A$ and any element $y\in A$, the equality
$x^* (y) = \1^* (y\odot x^{-1})$ holds, where $\1\in A$.
\end{proposition}

The mapping $A\times A\to \widehat K$ defined by the formula
$(x, y) \mapsto \langle x, y\rangle = \1^* (x\odot y)$
is called the  {\it canonical scalar product} (or simply {\it scalar product}).
The basic properties of the scalar product are easily derived from
Proposition~\ref{p:semialg} (in particular, the scalar product is commutative
if the $b$-semialgebra $A$ is commutative). The following theorem is an idempotent version of the Riesz--Fisher theorem.

\begin{theorem}
\label{t:RF}
Let a $b$-semialgebra $A$ be a semifield.
Then any nonzero $a$-linear functional $f$ on $A$ can be represented
as $f(y) = \langle y, x\rangle$, where $x\in A$, $x\neq \0$ and
$\langle \cdot, \cdot\rangle$ is the canonical scalar product on
$A$.
\end{theorem}

\begin{remark} Using the completion precedures, one can extend all
the results obtained to the case of incomplete semirings, spaces,
 and semimodules, see~\cite{LMS-01}.
\end{remark}

\begin{example} Let ${\mathcal B} (X)$ be a set of all bounded
functions with values belonging to $\R(\max, +)$ on an arbitrary
set $X$ and let $\widehat{\mathcal B}(X) = {\mathcal B} (X)\cup \{\0\}$.
The pointwise idempotent addition of functions
$(\varphi_1 \oplus \varphi_2) (x) = \varphi_1(x)\oplus \varphi_2(x)$
and the multiplication $(\varphi_1\odot \varphi_2)(x) = (\varphi_1(x))\odot
( \varphi_2(x))$ define on $\widehat{\mathcal B} (X)$ the structure of a
$b$-semialgebra over the $b$-complete semifield $\R_{\max}$. In
this case, $\1^* (\varphi) =\sup_{x\in X} \varphi (x)$
and the scalar product is expressed in terms of idempotent integration: $\langle\varphi_1, \varphi_2\rangle  =
\sup_{x\in X} (\varphi_1 (x)\odot\varphi_2 (x)) =
\sup_{x\in X} (\varphi_1 (x) + \varphi_2 (x))
=\int\limits^{\oplus}_X(\varphi_1 (x)\odot\varphi_2 (x))\; dx$.
Scalar products of this type were systematically used in
idempotent snslysis. Using Theorems~\ref{t:xfunc} and~\ref{t:RF},
one can easily
describe $a$-linear functionals on idempotent spaces in terms of
idempotent measures and integrals.
\end{example}

\begin{example} Let $X$ be a linear space in the traditional sense.
The idempotent semiring (and linear space over $\R(\max, +)$) of convex functions
Conv$(X, \R)$
is $b$-complete but it is not a $b$-semialgebra
over the
semifield
$K = \R(\max, +)$.

 Any nonzero $a$-linear functional~$f$
on
Conv$(X, \R)$
has the form
$$
{\varphi}\mapsto f({\varphi})
=\sup_x\{{\varphi}(x)+\psi(x)\}
=\int^\oplus_X{\varphi(x)}\odot\psi(x)\,dx,
$$
where~$\psi$
is a concave function, i.e., an element of the idempotent
space Conc($X$, $\R$) = - Conv($X$, $\R$).
\end{example}

\subsection{Linear operator, $b$-semimodules and subsemimodules}
\label{ss:linoper}

In what follows, we suppose that all semigroups, semirings, semifields, semimodules, and spaces are idempotent unless
otherwise specified. We fix a basic semiring $K$ and examine semimodules and subsemimodules over $K$. We suppose that every
linear functional takes it values in the basic semiring.

Let $V$ and $W$ be $b$-complete semimodules over a $b$-complete
semiring $K$. Denote by $L_b(V,W)$ the set of all $b$-linear
mappings from $V$ to $W$. It is easy to check that $L_b(V,W)$
is an idempotent semigroup with respect to the pointwise
addition of operators; the composition (product) of $b$-linear
operators is also a $b$-linear operator, and therefore
the set $L_b(V,V)$
is an idempotent semiring with respect to these operations,
see, e.g., \cite{LMS-01}. The following proposition can be treated
as a version of the Banach--Steinhaus theorem in idempotent
analysis (as well as Proposition~\ref{p:BS1} above).

\begin{proposition} Assume that $S$ is a subset in $L_b(V,W)$ and the set
$\{g(v)\mid g\in S\}$ is bounded in $W$ for every element
$v\in V$; thus the element $f(v)$ = $\sup_{g\in S}{g(v)}$
exists, because the semimodule $W$ is $b$-complete. Then the
mapping $v\mapsto f(v)$ is a $b$-linear operator, i.e., an
element of $L_b(V,W)$. The subset $S$ is bounded; moreover, $\sup S = f$.
\end{proposition}

\begin{corollary}
\label{c:bsemigroup}
The set $L_b(V,W)$ is a $b$-complete
idempotent semigroup
with respect to the (idempotent) pointwise addition of
operators. If $V = W$, then $L_b(V,V)$ is a $b$-complete idempotent
semiring with respect to the operations of pointwise addition and
composition of operators.
\end{corollary}

\begin{corollary} A subset $S$ is bounded in
$L_b(V,W)$  if and only if the set $\{g(v)\mid g\in S\}$ is
bounded in the semimodule $W$ for every element $v\in V$.
\end{corollary}

A subset of an idempotent semimodule is called a {\it subsemimodule}
if it is closed under addition and multiplication by scalar coefficients.
A subsemimodule $V$ of a $b$-complete semimodule $W$ is {\it b-closed}
if $V$ is closed under sums of any subsets of $V$
that are bounded in $W$.
A subsemimodule of a $b$-complete semimodule is called a
{\it b-subsemimodule} if the corresponding embedding is a $b$-homomorphism.
It is easy to see that each $b$-closed subsemimodule is a $b$-subsemimodule,
but the converse is not true.
The main feature of $b$-subsemimodules is that restrictions of $b$-linear
operators and functionals to these semimodules are $b$-linear.

{\it The following definitions are very important} for our
purposes. Assume that $W$ is an idempotent $b$-complete semimodule over a
$b$-complete idempotent semiring $K$ and $V$ is a subset of $W$ such
that $V$ is closed under multiplication by scalar coefficients
and is an upper semilattice with respect to the order induced
from $W$. Let us define an addition operation in $V$ by the
formula $x\oplus y = \sup\{ x, y \}$, where $\sup$ means the least
upper bound in $V$. If $K$ is a semifield, then $V$ is a
semimodule over $K$ with respect to this addition.

For an arbitrary $b$-complete semiring $K$, we will say that
$V$ is a {\it quasisubsemimodule} of $W$ if $V$ is a
semimodule with respect to this addition (this means that the
corresponding distribution laws hold).

Recall that the simbol $\wedge$ means the greatest lower bound (see Subsection 7.1 above). A quasisubsemimodule $V$ of an idempotent $b$-complete semimodule
$W$ is called a $\wedge$-{\it subsemimodule} if it contains
$\0$ and is closed under the operations of taking infima (greatest lower
bounds) in $W$. It is easy to check  that {\it each $\wedge$-subsemimodule is a $b$-complete
semimodule}.

Note that quasisubsemimodules and $\wedge$-subsemimodules
may fail to be subsemimodules, because only the order is induced
and not the
corresponding addition (see Example~\ref{ex:conc} below).

Recall that idempotent semimodules over semifields
are {\it idempotent spaces}. In idempotent mathematics, such spaces
are analogs of traditional linear (vector) spaces over fields.
In a similar way we use the corresponding terms like
{\it b-spaces, b-subspaces, b-closed subspaces},
$\wedge$-{\it subspaces}, etc.

Some examples are presented below.

\subsection{Functional semimodules}
\label{ss:funcsem}
Let $X$ be an arbitrary nonempty set and $K$ be an idempotent semiring.
By $K(X)$ denote the semimodule of all mappings (functions) $X \to K$
endowed with the pointwise operations. By $K_b(X)$ denote the subsemimodule
of $K(X)$ consisting of  all bounded
mappings. If $K$ is a $b$-complete semiring, then $K(X)$ and $K_b(X)$ are
$b$-complete semimodules. Note that $K_b(X)$ is a $b$-subsemimodule but not
a $b$-closed subsemimodule of $K(X)$. Given a point $x\in X$, by
$\delta _x$ denote the functional on $K(X)$ that maps
$f$ to $f(x)$. It can easily be checked that the functional
$\delta _x$ is $b$-linear on~$K(X)$.

Recall that the functional $\delta_x$ is generated by the usual function
\begin{equation*}
\delta_x(y)=
\begin{cases}
\1, & \text{if $x=y$},\\
\0, & \text{if $x\neq y$},
\end{cases}
\end{equation*}
so $\varphi(x)=\int^{\oplus}\delta_x(y)\varphi(y) dy=
\sup\limits_y (\delta_x(y)\odot\varphi(y))$. Note that $\delta$-functions
form a natural (continuous in general) basis in any typical
functional semimodule.

We say that a quasisubsemimodule of $K(X)$ is an (idempotent)
{\it functional semimodule} on the set $X$. An idempotent functional
semimodule in $K(X)$ is called {\it b-complete} if it is a $b$-complete
semimodule.

A functional semimodule $V\subset K(X)$ is called a {\it functional
b-semimodule} if it is a b-subsemimodule of $K(X)$; a functional
semimodule $V\subset K(X)$ is called a {\it functional $\wedge$-semimodule}
if it is a $\wedge$-subsemimodule of $K(X)$.

In general, a functional of the form $\delta _x$ on a functional semimodule
is not even linear, much less $b$-linear
(see Example~\ref{ex:conc} below). However, the following
proposition holds, which is a direct consequence of our definitions.

\begin{proposition} An arbitrary $b$-complete functional semimodule
$W$ on a set $X$ is a $b$-subsemimodule of $K(X)$ if and only if
each functional of the form $\delta _x$ (where $x\in X$) is
$b$-linear on $W$.
\end{proposition}

\begin{example}
\label{ex:bounded}
The semimodule $K_b(X)$ (consisting of all bounded
mappings from an arbitrary set $X$ to a
$b$-complete idempotent semiring $K$)
is a functional $\wedge$-semimodule. Hence it is a $b$-complete
semimodule over $K$. Moreover, $K_b(X)$ is a $b$-subsemimodule of the
semimodule $K(X)$ consisting of all mappings $X\to K$.
\end{example}

\begin{example}
\label{ex:finite}
If $X$ is a finite set consisting of $n$
elements ($n>0$), then $K_b(X) = K(X)$ is an
``$n$-dimensional'' semimodule over $K$; it is denoted by
$K^n$. In particular, $\R_{max}^n$ is an idempotent space over
the semifield $\R_{max}$, and $\widehat{\R}_{\max}^n$
is a semimodule over the semiring $\widehat{\R}_{\max}$. Note
that $\widehat{\R}_{\max}^n$ can be treated as a space
over the semifield $\R_{max}$. For example, the semiring
 $\widehat{\R}_{\max}$ can be treated as a space
(semimodule) over~$\Rmax$.
\end{example}

\begin{example}
\label{ex:USC}
Let $X$ be a topological space. Denote by
$USC(X)$ the set of all upper semicontinuous functions with
values in $\R_{\max}$. By definition, a function $f(x)$ is upper
semicontinuous if the set $X_s =\{ x\in X\mid f(x)\geq s\}$
is closed in $X$ for every element $s\in \R_{\max}$ (see, e.g.,
\cite{LMS-01}, Sec.~2.8). If a family $\{f_{\alpha}\}$ consists of
upper semicontinuous (e.g., continuous) functions and
$f(x) = \inf_{\alpha} f_{\alpha} (x)$, then $f(x) \in USC(X)$.
It is easy to check that $USC(X)$  has a natural structure
of an idempotent space over $\Rmax$. Moreover, $USC(X)$ is a functional
$\wedge$-space on $X$ and a b-space. The subspace
$USC(X)\cap K_b(X)$ of $USC(X)$ consisting of bounded (from
above)
functions has the same properties.
\end{example}

\begin{example}
\label{ex:conc}
 Note that an idempotent functional semimodule
(and even a functional $\wedge$-semimodule) on a set $X$ is not
necessarily a subsemimodule of $K(X)$. The simplest example is the
functional space (over $K=\Rmax$) Conc(${\R}$) consisting of all
concave functions on $\R$ with values in $\Rmax$.
Recall that a function $f$ belongs to Conc(${\R}$)
if and only if the subgraph
of this function is convex, i.e., the formula
$f(ax+(1-a)y)\geq af(x)+ (1-a)f(y)$ is valid for $0\leq a\leq 1$.
The basic operations with $\0\in \Rmax$ can be defined in an obvious way.
If $f,g \in$Conc$({\R})$, then denote by $f\oplus g$ the sum of
these functions in Conc$({\R})$. The subgraph of $f\oplus g$ is the
convex hull of the subgraphs of $f$ and $g$. Thus $f\oplus g$
does not coincide with the pointwise sum (i.e., $\max\{f(x), g(x)\}$).
\end{example}

\begin{example}
\label{ex:lip}
Let $X$ be a nonempty metric space with a
fixed metric $r$. Denote by Lip$(X)$ the set of all functions
defined on $X$ with values in $\Rmax$ satisfying the
following {\it Lipschitz condition}:
$$
\mid f(x)\odot (f(y))^{-1}\mid   =  \mid f(x) - f(y)\mid  \leq  r(x, y),
$$
where $x$, $y$ are arbitrary elements of $X$.  The set
Lip$(X)$ consists of continuous real-valued functions (but
not all of them!) and (by definition) the function equal to
$-\infty = \0$ at every point $x\in X$. The set Lip$(X)$ has the structure of an idempotent
space over the semifield $\Rmax$. Spaces of the form
Lip$(X)$ are said to be  {\it Lipschitz spaces}. These spaces are $b$-subsemimodules in $K(X)$.
\end{example}

\subsection{Integral representations of
linear operators in functional semimodules}

Let $W$ be an idempotent $b$-complete semimodule over a $b$-complete
semiring $K$ and $V\subset K(X)$ be a $b$-complete functional
semimodule on $X$. A mapping $A:V\to W$ is called an {\it integral operator}
or an operator with an {\it integral representation} if there exists
a mapping $k:X\to W$, called the {\it integral kernel} (or {\it kernel})
{\it of the operator} $A$, such that
\begin{equation}
\label{e:intker}
Af = \sup_{x\in X}{(f(x)\odot k(x))}.
\end{equation}
In idempotent analysis, the right-hand side of formula (11) is often
written as $\int_X^{\oplus}f(x)\odot k(x) dx$.
Regarding the kernel $k$, it is assumed that the set
$\{f(x)\odot k(x)|x\in X\}$ is bounded in $W$ for all $f\in V$ and $x\in X$.
We denote the set of all functions with this property by
${\text{\rm kern}}_{V,W}(X)$. In particular, if $W=K$ and $A$ is a functional,
then this functional is called {\it integral}. Thus each integral
functional can be presented in the form of a ``scalar product''
$f \mapsto \int_X^{\oplus} f(x) \odot k(x)\; dx$, where
$k(x)\in K(X)$; in idempotent analysis, this situation is
standard.

Note that a functional of the form $\delta_y$ (where $y\in X$)
is a typical integral functional; in this case,
$k(x) = \1$ if $x = y$ and $k(x) = \0$ otherwise.

We call a functional semimodule $V\subset K(X)$ {\it nondegenerate}
if for every point $x\in X$ there exists a function $g\in V$
such that $g(x)=\1$, and {\it admissible} if for every function $f\in V$
and every point $x\in X$ such that $f(x)\neq \0$ there exists a
function $g\in V$ such that $g(x)=\1$ and~$f(x)\odot g\preceq f$.

Note that all idempotent functional semimodules over semifields
are admissible (it is sufficient to set $g = f(x)^{-1}\odot f$).

\begin{proposition} Denote by $X_V$  the subset of $X$ defined by the
formula $X_V=\{ x\in X\mid\; \exists f\in V: f(x)=\1 \}$. If the semimodule
$V$ is admissible, then the restriction to $X_V$ defines an
embedding $i:V\to K(X_V)$ and its image $i(V)$ is admissible and
nondegenerate.

If a mapping $k:X\to W$ is a kernel of a mapping $A:V\to W$, then the
mapping $k_V:X\to W$ that is equal to $k$ on $X_V$ and equal to $\0$
on $X\smallsetminus~X_V$ is also a kernel of~$A$.

A mapping $A:V\to W$ is integral if and only if the mapping
 $i_{-1}A:i(A)\to W$ is integral.
\end{proposition}

In what follows, $K$ always denotes a fixed $b$-complete idempotent
(basic) semiring.  If an operator has an integral representation,
this representation may not be unique. However, if the semimodule
$V$ is nondegenerate, then the set of all kernels of a fixed
integral operator is bounded with respect to the natural order
in the set of all kernels and is closed under the supremum
operation applied to its arbitrary subsets. In particular,
{\it any integral
operator defined on a nondegenerate functional semimodule has
a unique maximal kernel}.

An important point is that an integral operator is not necessarily $b$-linear
and even linear except when $V$ is a $b$-subsemimodule of $K(X)$
(see Proposition~\ref{p:funcsem} below).

If $W$ is a functional semimodule on a nonempty set $Y$,
then an integral kernel $k$ of an operator $A$ can be naturally identified
with the function on $X\times Y$ defined by the formula
$k(x,y)=(k(x))(y)$. This function will also be called an
{\it integral kernel} (or {\it kernel}) of the operator $A$.
As a result, the set  ${\text{\rm kern}}_{V,W}(X)$  is identified with the set
${\text{\rm kern}}_{V,W}(X,Y)$ of all mappings $k: X\times Y\to K$ such that for every
point $x\in X$ the mapping $k_x: y\mapsto k(x,y)$ lies in $W$ and for
every $v\in V$ the set $\{v(x)\odot k_x|x\in X\}$ is bounded in $W$.
Accordingly, the set of all integral kernels of $b$-linear
operators can be embedded into~${\text{\rm kern}}_{V,W}(X,Y)$.

If $V$ and $W$ are functional $b$-semimodules on $X$ and $Y$,
respectively, then the set of all kernels of $b$-linear operators can
be identified with ${\text{\rm kern}}_{V,W}(X,Y)$
and the following formula holds:
\begin{equation}
\label{e:intrepr}
Af(y)=\sup_{x\in X}{(f(x)\odot k(x,y))}=\int_X^{\oplus} f(x)
\odot k(x,y) dx.
\end{equation}
This formula coincides with the usual definition of an
integral representation of an operator. Note that formula~\eqref{e:intker} can
be rewritten in the form
\begin{equation}
\label{e:intker2}
Af=\sup_{x\in X} {(\delta _x(f)\odot k(x))}.
\end{equation}

\begin{proposition}
\label{p:funcsem}
An arbitrary b-complete functional semimodule
$V$ on a nonempty set $X$ is a functional b-se\-mi\-mo\-du\-le on
$X$ (i.e., a b-sub\-semi\-mo\-du\-le of $K(X)$) if and only if
all integral operators defined on $V$ are b-linear.
\end{proposition}

The following notion (definition) is especially important
for our purposes. Let $V\subset K(X)$ be a $b$-complete functional
semimodule over a $b$-complete idempotent semiring $K$. We
say that the {\it kernel theorem} holds for the
semimodule $V$ if every $b$-linear mapping from $V$ into an
arbitrary $b$-complete semimodule over $K$ has an integral
representation.

\begin{theorem}
\label{t:intrepr}
Assume that  a b-complete
semimodule $W$ over a b-complete semiring $K$ and an
admissible functional $\wedge$-semimodule $V\subset K(X)$ are given.
Then every b-linear operator $A:V\to W$ has an integral
representation of the form~\eqref{e:intker}. In particular, if $W$ is
a functional b-semimodule on a set $Y$, then the operator
$A$ has an integral representation of the form~\eqref{e:intrepr}. Thus
for the semimodule $V$ the kernel theorem holds.
\end{theorem}

\begin{remark} Examples of admissible functional
$\wedge$-semimodules (and $\wedge$-spaces) appearing in
Theorem~\ref{t:intrepr}
are presented above, see, e.g., examples~\ref{ex:bounded}  --
\ref{ex:USC}. Thus for these functional
semimodules and spaces $V$ over $K$, the kernel theorem holds
and every $b$-linear  mapping $V$ into an arbitrary $b$-complete
semimodule $W$ over $K$ has an integral representation~\eqref{e:intrepr}.
Recall that every functional space over a $b$-complete semifield is
admissible, see above.
\end{remark}

\subsection{Nuclear operators and their integral representations}
\label{ss:nuclear}

Let us introduce some important definitions. Assume that
$V$ and $W$ are $b$-complete semimodules. A mapping $g: V\to
W$ is called {\it one-dimensional} (or a {\it mapping of
rank} 1) if it is of the form $v\mapsto \phi (v)\odot w$,
where $\phi$ is a $b$-linear functional on $V$ and $w\in W$.
A mapping $g$ is called {\it b-nuclear} if it is the sum (i.e., supremum)
of a bounded set of one-dimensional mappings. Since every
one-dimensional mapping is $b$-linear (because the functional
$\phi$ is $b$-linear), {\it every b-nuclear operator is b-linear}
(see Corollary~\ref{c:bsemigroup} above). Of course, $b$-nuclear mappings are closely
related to tensor products of idempotent semimodules, see~\cite{LMS-99}.

By $\phi\odot w$ we denote the one-dimensional operator
$v\mapsto\phi(v)\odot w$. In fact, this is an element of
the corresponding tensor product.

\begin{proposition} The composition (product) of a
b-nuclear and a
b-linear mapping or of a b-linear and a b-nuclear mapping
is a b-nuclear operator.
\end{proposition}

\begin{theorem} Assume that $W$ is a b-complete semimodule over a
b-complete semiring $K$ and $V\subset K(X)$ is a functional b-semimodule.
If every b-linear functional on $V$ is integral, then a b-linear
operator $A:V\to W$ has an integral representation if and only if
it is b-nuclear.
\end{theorem}

\subsection{The $b$-approximation property and $b$-nuclear
 semimodules and spaces}

We say that a $b$-complete semimodule $V$ has the {\it
b-approximation property} if the identity operator id:$V\to V$ is
$b$-nuclear (for a treatment of the approximation property for
locally convex spaces in the traditional functional
analysis, see \cite{Scha}).

Let $V$ be an arbitrary $b$-complete semimodule over a $b$-complete
idempotent semiring $K$. We call this semimodule a {\it b-nuclear
semimodule} if any $b$-linear mapping of $V$ to an arbitrary
$b$-complete semimodule $W$ over $K$ is a $b$-nuclear operator. Recall
that, in the traditional functional analysis, a locally convex space
is nuclear if and only if all continuous linear mappings of this
space to any Banach space are nuclear operators, see~\cite{Scha}.

\begin{proposition} Let $V$ be an arbitrary
b-complete semimodule over a b-complete semiring $K$. The
following statements are equivalent:
\begin{itemize}
\item[1] the semimodule $V$ has the b-approximation property;
\item[2] every b-linear mapping from $V$ to an arbitrary b-complete
semimodule $W$ over $K$ is b-nuclear;
\item[3] every b-linear mapping from an arbitrary b-complete
semimodule $W$ over $K$ to the semimodule $V$ is b-nuclear.
\end{itemize}
\end{proposition}

\begin{corollary} An arbitrary b-complete semimodule over
a b-complete semiring $K$ is b-nuclear if and only if this
semimodule has the b-approximation property.
\end{corollary}

Recall that, in the traditional functional analysis, any
nuclear space has the approximation property but the converse
is not true.

Concrete examples of $b$-nuclear spaces and semimodules are
described in Examples~\ref{ex:bounded},~\ref{ex:finite} and~\ref{ex:lip} (see above). Important $b$-nuclear spaces and semimodules (e.g., the so-called Lipschitz spaces and
semi-Lipschitz semimodules) are described in \cite{LSz-07b}.
In this paper there is a description of all functional $b$-semimodules
for which the kernel theorem holds (as semi-Lipschitz semimodules);
this result is due to G.~B.~Shpiz.

It is
easy to show that the idempotent spaces $USC(X)$ and
Conc($\R$)  (see Examples~\ref{ex:USC} and~\ref{ex:conc}) are not $b$-nuclear
(however, for these spaces the kernel theorem is true). The reason
is that these spaces are not functional $b$-spaces and the
corresponding $\delta$-functionals are not $b$-linear
(and even linear).

\subsection{Kernel theorems for functional $b$-semimodules}

Let $V\subset K(X)$ be a $b$-complete functional semimodule
over a $b$-complete semiring $K$. Recall that for $V$ the
{\it kernel theorem} holds if every $b$-linear mapping of
this semimodule to an arbitrary $b$-complete semimodule over
$K$ has an integral representation.

\begin{theorem} Assume that  a b-complete
semiring $K$ and a nonempty set $X$ are given. The kernel theorem
holds for any functional b-semimodule $V\subset K(X)$ if
and only if every b-linear functional on $V$ is integral
and the semimodule $V$ is b-nuclear, i.e., has the
b-approximation property.
\end{theorem}

\begin{corollary} If for a functional b-semimodule
the kernel theorem holds, then this semimodule is b-nuclear.
\end{corollary}

Note that the possibility to obtain an integral representation of
a functional means that one can decompose it into a sum
of functionals of the form $\delta _x$.

\begin{corollary} Assume that  a b-complete semiring $K$
and a nonempty set $X$ are given. The kernel theorem holds for a
functional b-semimodule
$V\subset K(X)$ if and only if the identity operator
id: $V\to V$ is integral.
\end{corollary}

\subsection{Integral representations of operators in
 abstract idempotent semimodules}
\label{ss:intrepr}

In this subsection, we examine the following problem: when
a $b$-complete idempotent semimodule $V$ over a $b$-complete
semiring is isomorphic to a functional $b$-semimodule
$W$ such that the kernel theorem holds for $W$.

Assume that $V$ is a $b$-complete idempotent
semimodule over a $b$-complete semiring $K$ and $\phi$ is a
$b$-linear functional defined on $V$. We call this
functional a $\delta$-{\it functional} if there exists an
element $v\in V$ such that
$$
\phi(w)\odot v\preceq w
$$
for every element $w\in V$. It is easy to see that every
functional of the form $\delta_x$ is a $\delta$-functional
in this sense (but the converse is not true in general).

Denote by $\Delta(V)$ the set of all $\delta$-functionals on
$V$. Denote by $i_\Delta$ the natural mapping
$V\to K(\Delta(V))$ defined by the formula
$$
(i_\Delta (v) )(\phi)=\phi(v)
$$
for all $\phi \in \Delta(V)$.
We say that an element $v\in V$ is {\it pointlike} if there
exists a $b$-linear functional $\phi$ such that
$\phi(w)\odot v\preceq w$ for all $w\in V$. The set of all
pointlike elements of $V$ will be denoted by $P(V)$. Recall
that by $\phi\odot v$ we denote the one-dimensional operator
$w\mapsto \phi(w)\odot v$.

The following assertion is an
obvious consequence of our definitions (including the
definition of the standard order) and the idempotency
of our addition.

\begin{remark} If a one-dimensional operator $\phi\odot v$
appears in the decomposition of the identity operator on $V$
into a sum of one-dimensional operators, then $\phi\in\Delta(V)$
and $v\in P(V)$.
\end{remark}

Denote by $id$ and $Id$ the identity operators on $V$ and
$i_\Delta (V)$, respectively.

\begin{proposition} If the operator id is
b-nuclear, then $i_\Delta$ is an embedding and the
operator Id is integral.

If the operator  $i_\Delta$ is an embedding and
the operator Id is integral, then the operator id is
$b$-nuclear.
\end{proposition}

\begin{theorem} A b-complete idempotent semimodule $V$ over
a b-complete idempotent semiring $K$ is isomorphic to a
functional b-semimodule for which the kernel theorem holds
if and only if the identity mapping on $V$ is a
b-nuclear operator, i.e., $V$ is a b-nuclear semimodule.
\end{theorem}

The following proposition shows that, in a certain sense, the
embedding $i_\Delta$ is a universal representation of a
$b$-nuclear semimodule in the form of a functional
$b$-semimodule for which the kernel theorem holds.

\begin{proposition} Let $K$ be a b-complete
idempotent semiring, $X$ be a nonempty set, and $V\subset
K(X)$ be a functional $b$-semimodule on $X$ for which the
kernel theorem holds. Then there exists a natural mapping
$i:X\to \Delta(V)$ such that the corresponding mapping
$i_*: K(\Delta(V))\to K(X)$ is an isomorphism of
$i_\Delta(V)$ onto~$V$.
\end{proposition}

\section{The dequantization transform, convex geometry and the Newton polytopes}
\label{s:dequantfunc}

Let $X$ be a topological space. For functions $f(x)$ defined on
$X$ we shall say that a certain property is valid {\it almost
everywhere} (a.e.) if it is valid for all elements $x$ of an open
dense subset of $X$. Suppose $X$ is $\C^n$ or $\R^n$; denote by
$\R^n_+$ the set $x=\{\,(x_1, \dots, x_n)\in X \mid x_i\geq 0$ for
$i = 1, 2, \dots, n$.
 For $x= (x_1, \dots, x_n) \in X$ we set
${\mbox{exp}}(x) = ({\mbox{exp}}(x_1), \dots, {\mbox{exp}}(x_n))$;
so if $x\in\R^n$, then ${\mbox{exp}}(x)\in \R^n_+$.

Denote by $\cF(\C^n)$ the set of all functions defined and
continuous on an open dense subset $U\subset \C^n$ such that
$U\supset \R^n_+$. It is clear that $\cF(\C^n)$ is a ring (and an
algebra over $\C$) with respect to the usual addition and
multiplications of functions.

For $f\in \cF(\C^n)$ let us define the function $\hat f_h$ by the
following formula:
\begin{equation}
\label{e:logtrans}
\hat f_h(x) = h \log|f({\mbox{exp}}(x/h))|,
\end{equation}
where $h$ is a (small) real positive parameter and $x\in\R^n$. Set
\begin{equation}
\label{e:hatfx}
\hat f(x) = \lim_{h\to +0} \hat f_h (x),
\end{equation}
if the right-hand side of~\eqref{e:hatfx} exists almost everywhere.

We shall say that the function $\hat f(x)$ is a {\it
dequantization} of the function $f(x)$ and the map $f(x)\mapsto
\hat f(x)$ is a {\it dequantization transform}. By construction,
$\hat f_h(x)$ and $\hat f(x)$ can be treated as functions taking
their values in $\R_{\max}$. Note that in fact $\hat f_h(x)$ and
$\hat f(x)$
 depend on the restriction of $f$ to $\R_+^n $
only; so in fact the dequantization transform is constructed for
functions defined on $\R^n_+$ only. It is clear that the
dequantization transform is generated by the Maslov dequantization
and the map $x\mapsto |x|$.

Of course, similar definitions can be given for functions defined
on $\R^n$ and $\R_+^n$. If $s=1/h$, then we have the following
version of~\eqref{e:logtrans} and~\eqref{e:hatfx}:
\begin{equation}
\label{e:hatfx1}
\hat f(x) = \lim_{s\to \infty} (1/s) \log|f(e^{sx})|.
\end{equation}

Denote by $\partial \hat f$ the subdifferential of the function
$\hat f$ at the origin.

If $f$ is a polynomial  we have
\begin{equation*}
\partial \hat f = \{\, v\in \R^n\mid (v, x) \leq \hat f(x)\
\forall x\in \R^n\,\}.
\end{equation*}

It is well known that all the convex compact subsets in $\R^n$
form an idempotent semiring $\mathcal{S}$ with respect to the
Minkowski operations: for $\alpha, \beta \in \mathcal{S}$ the sum
$\alpha\oplus \beta$ is the convex hull of the union $\alpha\cup
\beta$; the product $\alpha\odot \beta$ is defined in the
following way: $\alpha\odot \beta = \{\, x\mid x = a+b$, where
$a\in \alpha, b\in \beta$, see Fig~\ref{f:fig3}. In fact $\mathcal{S}$ is an
idempotent linear space over $\R_{\max}$.

Of course, the Newton polytopes of polynomials in $n$ variables
form a subsemiring
$\mathcal{N}$ in $\mathcal{S}$. If $f$, $g$ are polynomials, then
$\partial(\widehat{fg}) = \partial\hat f\odot\partial\widehat g$;
moreover, if $f$ and $g$ are ``in general position'', then
$\partial(\widehat{f+g}) = \partial\hat f\oplus\partial\widehat
g$. For the semiring of all polynomials with nonnegative
coefficients the dequantization transform is a homomorphism of
this ``traditional'' semiring to the idempotent semiring
$\mathcal{N}$.

\begin{figure}[b]
\includegraphics[scale=1]{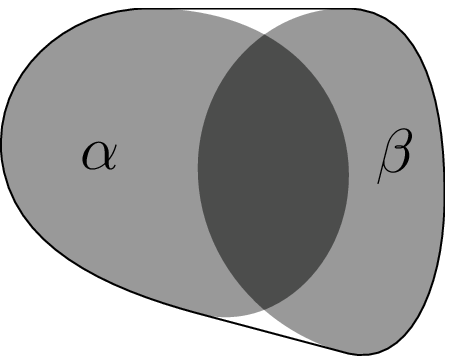}
%
%
\caption{Algebra of convex subsets.}
\label{f:fig3}       
\end{figure}

\if{
\begin{figure}
\centering
\epsfig{file=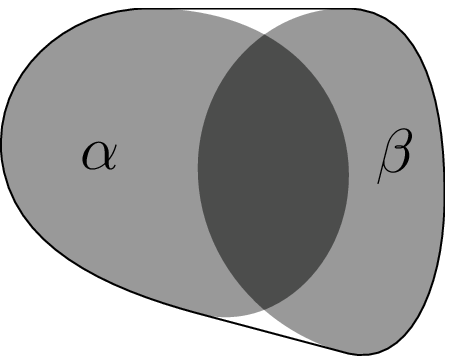,width=0.4\linewidth}
\caption{Algebra of convex subsets.}
\label{f:fig3}
\end{figure}
}\fi

\begin{theorem}
\label{t:newtpoly}
If $f$ is a polynomial, then the subdifferential $\partial\hat f$
of $\hat f$ at the origin coincides with the Newton polytope of
$f$.  For the semiring of polynomials with nonnegative
coefficients, the transform $f\mapsto\partial\hat f$ is a
homomorphism of this semiring to the semiring of convex polytopes
with respect to the Minkowski operations (see above).
\end{theorem}

Using the dequantization transform it is possible to generalize
this result to a wide class of functions and convex
sets, see below and~\cite{LSz-05}.

\subsection{Dequantization transform: algebraic properties}
\label{ss:dequantfunc}

Denote by $V$ the set $\rset^n$ treated as a linear Euclidean
space (with the scalar product $(x, y) = x_1y_1+ x_2y_2 +\dots + x_ny_n$)
and set $V_+ = \rset_+^n$.
We shall say that a function $f\in \maF(\cset^n)$ is {\it dequantizable}
whenever its
dequantization $\hat f(x)$ exists (and is defined on an open dense subset
of $V$). By $\maD (\cset^n)$ denote the set of all dequantizable functions and by
$\widehat{\maD}(V)$ denote the set $\{\,\hat f \mid f\in \maD(\cset^n)\,\}$. Recall that
functions from $\maD(\cset^n)$ (and $\widehat{\maD}(V)$) are defined almost everywhere and
$f=g$ means that $f(x) = g(x)$ a.e., i.e., for $x$ ranging over an open dense subset
of $\cset^n$ (resp., of $V$). Denote by $\maD_+(\cset^n)$ the set of all
functions $f\in \maD(\cset^n)$
such that $f(x_1, \dots, x_n)\geq 0$ if $x_i\geq 0$ for $i= 1,\dots, n$; so
 $f\in \maD_+(\cset^n)$ if the restriction of $f$ to $V_+ = \rset_+^n$ is a
 nonnegative function. By $\widehat{\maD}_+(V)$ denote the image of $\maD_+(\cset^n)$
 under the dequantization transform. We shall say that functions
$f, g\in \maD(\cset^n)$ are in {\it general position} whenever
$\hat f (x) \neq \widehat g(x)$ for $x$ running an open dense
subset of $V$.

\begin{theorem}
\label{t:dequantfunc}
For functions $f, g \in \maD(\cset^n)$ and any nonzero constant $c$, the
following equations are valid:

\begin{enumerate}
\item[1)] $\widehat{fg} = \hat f + \widehat g$;
\item[2)] $|\hat f| = \hat f$; $\widehat{cf} = f$; $\widehat c =0$;
\item[3)] $(\widehat{f+g})(x) = \max\{\hat f(x), \widehat g(x)\}$ a.e.\ if $f$ and $g$
are nonnegative on $V_+$ (i.e., $f, g \in \maD_+(\cset^n)$) or $f$ and $g$
 are in general position.
\end{enumerate}
Left-hand sides of these equations are well-defined automatically.
\end{theorem}

\begin{corollary}
\label{c:d+cn}
The set $\maD_+(\cset^n)$ has a natural structure of a semiring with respect to
the
 usual addition and multiplication of functions taking their values in $\cset$.
 The set $\widehat{\maD}_+(V)$ has a natural
 structure of an idempotent semiring with respect to the operations
$(f\oplus g)(x) = \max \{ f(x), g(x)\}$, $(f\odot g)(x) = f(x) + g(x)$;
elements
 of $\widehat{\maD}_+(V)$ can be naturally treated as functions taking their values
 in $\rset_{\max}$. The dequantization transform generates a homomorphism from
 $\maD_+(\cset^n)$ to $\widehat{\maD}_+(V)$.
\end{corollary}

\subsection{Generalized polynomials and simple functions}
\label{ss:genpoly}

For any nonzero number $a\in\cset$ and any vector
 $d = (d_1, \dots, d_n)\in V = \rset^n$
we set $m_{a,d}(x) = a \prod_{i=1}^n x_i^{d_i}$; functions of this kind we
 shall
call {\it generalized monomials}. Generalized monomials are defined a.e.\ on
$\cset^n$ and on $V_+$, but not on $V$ unless the numbers $d_i$ take
integer or suitable rational values. We shall say that a function $f$ is a {\it
generalized polynomial} whenever it is a finite sum of linearly
 independent generalized monomials. For instance, Laurent polynomials
and Puiseax polynomials are examples of generalized polynomials.

As usual, for $x, y\in V$ we set $(x,y) = x_1y_1 + \dots + x_ny_n$. The
following proposition is a result of a trivial calculation.

\begin{proposition}
\label{p:madhx}
For any nonzero number $a\in V = \cset$ and any vector $d\in V = \rset^n$
we have $(\widehat{m_{a,d}})_h(x) = (d, x) + h\log|a|$.
\end{proposition}

\begin{corollary}
\label{c:fhatlin}
If $f$ is a generalized monomial, then $\hat f$ is a linear function.
\end{corollary}

Recall that a real function $p$ defined on $V = \rset^n$ is {\it sublinear}
if $p = \sup_{\alpha}
p_{\alpha}$, where $\{p_{\alpha}\}$ is a collection of linear functions.
Sublinear functions defined everywhere on $V=\rset^n$ are convex; thus these
 functions are continuous, see \cite{MT:03}.
 We discuss sublinear functions of this kind only. Suppose $p$ is a continuous
 function defined on $V$, then $p$ is sublinear whenever

1) $p(x+ y) \leq p(x) + p(y)$ for all $x, y \in V$;

2) $p(cx) = cp(x)$ for all $x\in V$, $c\in \rset_+$.

So if $p_1$, $p_2$ are sublinear functions, then $p_1 +p_2$ is a sublinear
 function.

We shall say
that a function $f \in \maF(\cset^n)$ is {\it simple}, if its dequantization $\hat f$
exists and a.e.\ coincides with a sublinear function; by misuse of language, we
shall denote this (uniquely defined everywhere on $V$) sublinear
 function by the same symbol $\hat f$.

 Recall that simple functions $f$ and $g$ are {\it in general position} if
$\hat f(x) \neq \widehat g(x)$ for all $x$ belonging to an open dense subset of $V$.
In particular, generalized monomials are in
general position whenever they are linearly independent.

Denote by $\mathit{Sim}(\cset^n)$ the set of all simple functions defined on $V$ and
denote by $\mathit{Sim}_+(\cset^n)$ the set $\mathit{Sim}(\cset^n) \cap \maD_+(\cset^n)$. By
$\mathit{Sbl}(V)$ denote the
 set of all (continuous) sublinear functions defined on $V = \rset^n$ and by
 $\mathit{Sbl}_+(V)$ denote the image $\widehat{\mathit{Sim}_+}(\cset^n)$ of $\mathit{Sim}_+(\cset^n)$ under the
 dequantization transform.

The following statements can be easily deduced from Theorem 8.2 and definitions.

\begin{corollary}
The set $\mathit{Sim}_+(\cset^n)$ is a subsemiring of $\maD_+(\cset^n)$ and $\mathit{Sbl}_+(V)$
is an idempotent subsemiring of $\widehat{\maD_+}(V)$. The
dequantization transform generates an epimorphism of
$\mathit{Sim}_+(\cset^n)$ onto $\mathit{Sbl}_+(V)$. The set $\mathit{Sbl}(V)$ is an idempotent
semiring with respect to the operations
$(f\oplus g)(x) = \max \{ f(x), g(x)\}$,
$(f\odot g)(x) = f(x) + g(x)$.
\end{corollary}

\begin{corollary}
Polynomials and generalized polynomials are simple functions.
\end{corollary}

We shall say that functions $f, g\in\maD(V)$ are {\it asymptotically equivalent}
whenever $\hat f = \widehat g$; any simple function $f$ is an {\it asymptotic
monomial} whenever
$\hat f$ is a linear function. A simple function $f$ will be called an {\it
asymptotic polynomial} whenever $\hat f$ is a sum of a finite collection of
nonequivalent asymptotic monomials.

\begin{corollary}
Every asymptotic polynomial is a simple function.
\end{corollary}

\begin{example} Generalized polynomials, logarithmic functions of
(generalized) polynomials, and products of
polynomials and logarithmic functions are asymptotic polynomials. This follows
from our definitions and formula~\eqref{e:hatfx}.
\end{example}
\medskip

\subsection{Subdifferentials of sublinear functions}
\label{ss:subdiff}
We shall use some elementary results from convex analysis. These results can be
found, e.g., in \cite{MT:03}, ch. 1, \S 1.

For any function $p\in \mathit{Sbl}(V)$ we set
\begin{equation}
\partial p = \{\, v\in V\mid (v, x) \leq p(x)\ \forall x\in V\,\}.
\end{equation}

It is well known from convex analysis that for any sublinear function $p$ the
set $\partial p$ is exactly the {\it subdifferential} of $p$ at the origin.
 The following propositions are also known in convex
analysis.

\begin{proposition}
\label{p:subdiff}
Suppose $p_1,p_2\in \mathit{Sbl}(V)$, then
\begin{enumerate}
\item[1)] $\partial (p_1+p_2) = \partial p_1\odot\partial p_2 = \{\, v\in V\mid
v = v_1+v_2, \text{ where $v_1\in \partial p_1, v_2\in \partial p_2$}\,\}$;
\item[2)] $\partial (\max\{p_1(x), p_2(x)\}) = \partial p_1\oplus\partial p_2$.
\end{enumerate}
\end{proposition}

Recall that $\partial p_1\oplus \partial p_2$ is a convex hull of the set
$\partial p_1\cup \partial p_2$.

\begin{proposition}
\label{p:partialp}
Suppose $p\in \mathit{Sbl}(V)$.  Then $\partial p$ is a nonempty convex compact
subset of $V$.
\end{proposition}

\begin{corollary}
\label{p:homom}
The map $p\mapsto \partial p$ is a homomorphism of the idempotent semiring
 $\mathit{Sbl}(V)$ (see Corollary~\ref{c:d+cn}) to the idempotent semiring $\mathcal{S}$ of all convex
 compact subsets of $V$ (see Subsection~\ref{ss:dequantfunc} above).
\end{corollary}

\subsection{Newton sets for simple functions}

For any simple function $f\in \mathit{Sim}(\cset^n)$ let us denote by $N(f)$ the set
 $\partial(\hat f)$. We shall call $N(f)$ the {\it Newton set} of the
 function $f$.

\begin{proposition}
For any simple function $f$, its Newton set $N(f)$ is a nonempty convex
 compact subset of $V$.
\end{proposition}

This proposition follows from Proposition~\ref{p:partialp} and definitions.

\begin{theorem}
\label{t:simpfunc}
Suppose that $f$ and $g$ are simple functions. Then
\begin{enumerate}
\item[1)] $N(fg) = N(f)\odot N(g) = \{\, v\in V\mid v = v_1 +v_2
 \text{ with
 $v_1 \in N(f), v_2 \in N(g)$}\,\}$;
\item[2)] $N(f+g) = N(f)\oplus N(g)$, if $f_1$ and $f_2$ are in general
 position or $f_1, f_2 \in \mathit{Sim}_+(\cset^n)$ {\rm (}recall that $N(f)\oplus N(g)$
is the convex hull of $N(f)\cup N(g)${\rm )}.
\end{enumerate}
\end{theorem}

This theorem follows from Theorem~\ref{t:dequantfunc},
Proposition~\ref{p:subdiff} and definitions.

\begin{corollary}
The map $f\mapsto N(f)$ generates a homomorphism from $\mathit{Sim}_+(\cset^n)$ to
$\mathcal{S}$.
\end{corollary}

\begin{proposition}
\label{p:prod}
Let $f = m_{a,d}(x) = a \prod^n_{i=1} x_i^{d_i}$ be a monomial; here
$d = (d_1, \dots, d_n) \in V= \rset^n$ and $a$ is a nonzero
complex number. Then $N(f) = \{ d\}$.
\end{proposition}

This follows from Proposition~\ref{p:madhx}, Corollary~\ref{c:fhatlin} and definitions.

\begin{corollary}
Let $f = \sum_{d\in D} m_{a_d,d}$ be a polynomial. Then $N(f)$ is the polytope
$\oplus_{d\in D}\{d\}$, i.e.\ the convex hull of the finite set $D$.
\end{corollary}

This statement follows from Theorem~\ref{t:simpfunc} and
Proposition~\ref{p:prod}. Thus in this case
 $N(f)$ is the well-known classical Newton polytope of the polynomial $f$.

Now the following corollary is obvious.

\begin{corollary}
Let $f$ be a generalized or asymptotic polynomial. Then its Newton set
 $N(f)$ is a convex polytope.
\end{corollary}

\begin{example}. Consider the one dimensional case, i.e., $V = \rset$ and
suppose
 $f_1 = a_nx^n + a_{n-1}x^{n-1} + \dots + a_0$ and $f_2 = b_mx^m + b_{m-1}
 x^{m-1} + \dots + b_0$, where $a_n\neq 0$, $b_m\neq 0$, $a_0 \neq 0$,
 $b_0 \neq 0$. Then $N(f_1)$ is the segment $[0, n]$ and $N(f_2)$ is the
 segment $[0, m]$. So the map $f\mapsto N(f)$ corresponds to the map
 $f\mapsto \deg (f)$, where $\deg(f)$ is a degree of the polynomial $f$. In
 this case Theorem 2 means that $\deg(fg) = \deg f + \deg g$ and
 $\deg (f+g) = \max \{\deg f, \deg g\} = \max \{n, m\}$ if $a_i\geq 0$,
 $b_i\geq 0$ or $f$ and $g$ are in general position.
\end{example}
\medskip

\section{Dequantization of set functions and measures on metric spaces}
\label{s:meas}

The following results are presented in~\cite{LSz-07a}.

\begin{example} Let $M$ be a metric space, $S$ its arbitrary
subset with a compact closure. It is well-known that a Euclidean
$d$-dimensional ball $B_{\rho}$ of radius $\rho$ has volume
$$
\ovol\nolimits_d(B_{\rho})=\frac{\Gamma(1/2)^d}{\Gamma(1+d/2)}\rho^d,
$$
where $d$ is a natural parameter. By means of this formula it is
possible to define a volume of $B_{\rho}$ for any {\it real} $d$.
Cover $S$ by a finite number of balls of radii $\rho_m$. Set
$$
v_d(S):=\lim_{\rho\to 0} \inf_{\rho_m<\rho} \sum_m
\ovol\nolimits_d(B_{\rho_m}).
$$
Then there exists a number $D$ such that $v_d(S)=0$ for $d>D$ and
$v_d(S)=\infty$ for $d<D$. This number $D$ is called the {\it
Hausdorff-Besicovich dimension} (or {\it HB-dimension}) of $S$,
see, e.g.,~\cite{Mas-07}. Note that a set of non-integral
HB-dimension is called a fractal in the sense of B.~Mandelbrot.
\end{example}

\begin{theorem}
Denote by $\cN_{\rho}(S)$ the minimal number of balls of radius
$\rho$ covering $S$. Then
$$
D(S)=\mathop{\underline{\lim}}\limits_{\rho\to +0} \log_{\rho}
(\cN_{\rho}(S)^{-1}),
$$
where $D(S)$ is the HB-dimension of $S$. Set $\rho=e^{-s}$, then
$$
D(S)=\mathop{\underline{\lim}}\limits_{s\to +\infty} (1/s) \cdot
\log \cN_{exp(-s)}(S).
$$
So the HB-dimension $D(S)$ can be treated as a result of a
dequantization of the set function $\cN_{\rho}(S)$.
\end{theorem}

\begin{example} Let $\mu$ be a set function on $M$ (e.g., a
probability measure) and suppose that $\mu(B_{\rho})<\infty$ for
every ball $B_{\rho}$. Let $B_{x,\rho}$ be a ball of radius $\rho$
having the point $x\in M$ as its center. Then define
$\mu_x(\rho):=\mu(B_{x,\rho})$ and let $\rho=e^{-s}$ and
$$
D_{x,\mu}:=\mathop{\underline{\lim}}\limits_{s\to +\infty}
-(1/s)\cdot\log (|\mu_x(e^{-s})|).
$$
This number could be treated as a dimension of $M$ at the point
$x$ with respect to the set function $\mu$. So this dimension is a
result of a dequantization of the function $\mu_x(\rho)$, where
$x$ is fixed. There are many dequantization procedures of this
type in different mathematical areas. In particular, V.P.~Maslov's
negative dimension (see~\cite{Mas-07}) can be treated similarly.
\end{example}
\medskip

\section{Dequantization of geometry}
\label{s:geom}

An idempotent version of real algebraic geometry was discovered in
the report of O.~Viro for the Barcelona Congress~\cite{Vir-00}.
Starting from the idempotent correspondence principle O.~Viro
constructed a piecewise-linear geometry of polyhedra of a special
kind in finite dimensional Euclidean spaces as a result of the
Maslov dequantization of real algebraic geometry. He indicated
important applications in real algebraic geometry (e.g.,
 in the framework of Hilbert's 16th problem for
constructing real algebraic varieties with prescribed properties
and parameters) and relations to complex algebraic geometry and
amoebas in the sense of I.~M.~Gelfand, M.~M.~Kapranov, and
A.~V.~Zelevinsky, see~\cite{GKZ,Vir-02}. Then complex
algebraic geometry was dequantized by G.~Mikhalkin and the result
turned out to be the same; this new `idempotent' (or asymptotic)
geometry is now often called the {\it tropical algebraic
geometry}, see,
e.g.,~\cite{IMS:07,LM:05,LMS:07,LS:09,Mik-05,Mik-06}.

There is a natural relation between the Maslov dequantization and
amoebas.

Suppose $({\cset}^*)^n$ is a complex torus, where ${\cset}^* =
{\cset}\backslash \{0\}$ is the group of nonzero complex numbers
under multiplication.  For
 $z = (z_1, \dots, z_n)\in
(\cset^*)^n$ and a positive real number $h$ denote by $\Log_h(z) =
h\log(|z|)$ the element
\[(h\log |z_1|, h\log |z_2|, \dots,
h\log|z_n|) \in \rset^n.\] Suppose $V\subset (\cset^*)^n$ is a
complex algebraic variety; denote by $\maA_h(V)$ the set
$\Log_h(V)$. If $h=1$, then the set $\maA(V) = \maA_1(V)$ is
called the {\it amoeba} of $V$; the amoeba $\maA(V)$ is a closed
subset of $\rset^n$ with a non-empty complement. Note that this
construction depends on our coordinate system.

For the sake of simplicity suppose $V$ is a hypersurface
in~$(\cset^*)^n$ defined by a polynomial~$f$; then there is a
deformation $h\mapsto f_h$ of this polynomial generated by the
Maslov dequantization and $f_h = f$ for $h = 1$. Let $V_h\subset
({\cset}^*)^n$ be the zero set of $f_h$ and set $\maA_h (V_h) =
{\Log}_h (V_h)$. Then
 there exists a tropical variety
$\mathit{Tro}(V)$ such that the subsets $\maA_h(V_h)\subset
\rset^n$ tend to $\mathit{Tro}(V)$ in the Hausdorff metric as
$h\to 0$. The tropical variety $\mathit{Tro}(V)$ is a result of a
deformation of the amoeba $\maA(V)$ and the Maslov dequantization
of the variety $V$. The set $\mathit{Tro}(V)$ is called the {\it
skeleton} of $\maA(V)$.

\begin{figure}[b]
\includegraphics[scale=1]{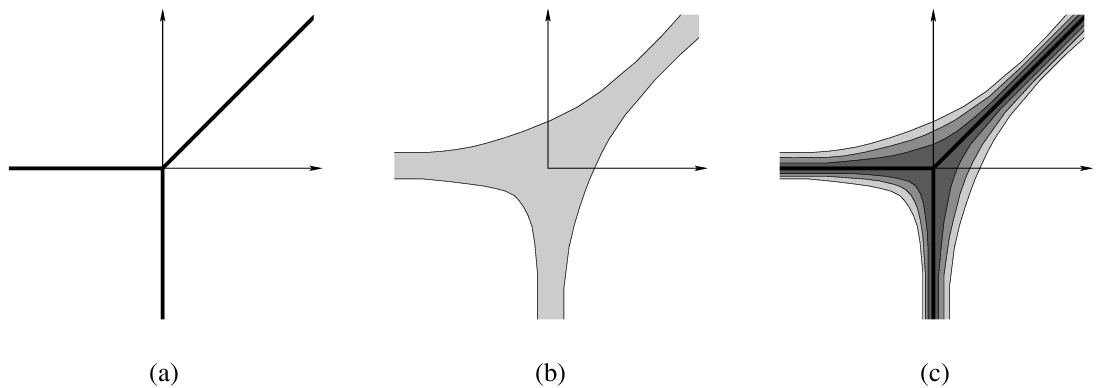}
%
%
\caption{Tropical line and deformations of an amoeba}
\label{f:amoeba}       
\end{figure}

\if{
\begin{figure}
\noindent\epsfig{file=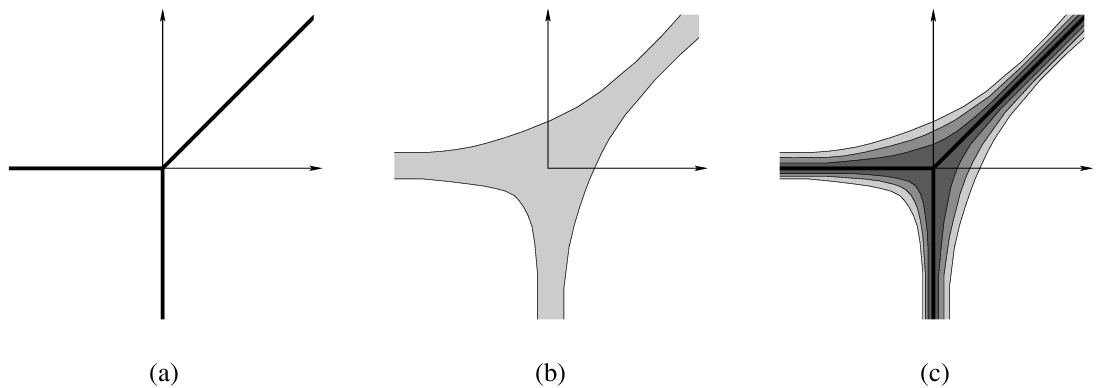,width=0.9\linewidth}
\caption{Tropical line and deformations of an amoeba.}
\end{figure}
}\fi

\begin{example}  For the line $V = \{\, (x, y)\in ({\cset}^*)^2
\mid x + y + 1 = 0\,\}$ the piecewise-linear graph
$\mathit{Tro}(V)$ is a tropical line, see Fig.~\ref{f:amoeba}(a). The amoeba
$\maA(V)$ is represented in Fig.~\ref{f:amoeba}(b), while Fig.~\ref{f:amoeba}(c) demonstrates
the corresponding deformation of the amoeba.
\end{example}

\section{Some semiring constructions and the matrix Bellman equation}
\label{s:bellman}

\subsection{Complete idempotent semirings and examples}

Recall that a partially ordered set $S$ is {\em complete} if for every subset
$T\subset S$ there exist elements $\sup T\in S$ and $\inf T\in S$. We say that an idempotent semiring $S$ is {\em complete} if it is complete as an ordered set with respect to
the standard order. Of course, any a-complete semiring (see subsect.~\ref{ss:semimod}) is complete. The most well-known and important examples are ``numerical
semirings'' consisting of (a subset of) real numbers and ordered by the usual linear order $\leq$.

\begin{example}
\label{ex:rmaxhat}
Consider the semiring $\rmaxh = \rmax \cup \{\infty\}$
with standard operations $\oplus=\max$, $\odot=+$ and
neutral elements $\0=-\infty$, $\1=0$, $x\leq\infty$,
$x \oplus \infty = \infty$ for all $x$, $x \odot
\infty = \infty \odot x = \infty$ if $x \neq \0$, and $\0 \odot \infty =
\infty \odot \0$. The semiring $\rmaxh$ is complete and a-complete. The semiring
$\rminh=\rmin\cup\{-\infty\}$ with obvious operations is also complete; $\rminh$ and $\rmaxh$ are isomorphic.
\end{example}

\begin{example}
\label{ex:maxmin}
Consider the semiring $S_{\max,\min}^{[a,b]}$ defined on the real interval
$[a,b]$ with operations $\oplus=\max$, $\odot=\min$ and neutral elements
$\0=a$ and $\1=b$. The semiring is complete and a-complete. Set $S_{\max,\min}=S_{\max,\min}^{[a,b]}$ with $a=-\infty$ and $b=+\infty$.
If $-\infty\leq a<b\leq +\infty$ then $S_{\max,\min}^{[a,b]}$  and
$S_{\max,\min}$ are isomorphic.
\end{example}

\begin{example}
\label{ex:boolean}
The Boolean algebra $B=\{\0,\1\}$ is a complete and a-complete semifield consisting of two elements.
\end{example}

\subsection{Closure operations}
\label{ss:closure}

Let a semiring~$S$ be endowed with a partial unary \emph{closure (or Kleene)
operation}~$*$ such that $x \preceq y$ implies $x^* \preceq y^*$ and $x^* = \1
\oplus (x^* \odot x) = \1 \oplus (x \odot x^*)$ on its domain of
definition. In particular, $\0^* = \1$ by definition. These axioms imply
that $x^* = \1 \oplus x \oplus x^2 \oplus \dots \oplus (x^* \odot x^n)$ if
$n \geqslant 1$. Thus $x^*$ can be considered as a `regularized sum' of the
series $x^* = \1 \oplus x \oplus x^2 \oplus \dots$; in an idempotent
semiring, by definition, $x^* = \sup \{ \1, x, x^2, \dots \}$ if this
supremum exists. So if $S$ is complete, then the closure operation is well-defined for every element $x\in S$.

In numerical semirings the operation~$*$ is defined as follows:
$x^* = (1-x)^{-1}$ if $x \prec 1$ in $\rset_+$, or $\widehat{\textbf{R}}_+$ and $x^*=\infty$ if $x\succcurlyeq 1$ in $\widehat{\textbf{R}}_+$; $x^* = \1$ if $x \preceq \1$ in $\rmax$
and $\rmaxh$, $x^* = \infty$ if $x \succ \1$ in $\rmaxh$, $x^* = \1$
for all $x$ in $\smaxmin^{[a,b]}$. In all other cases $x^*$ is undefined.
Note that the closure operation is very easy to implement.

\subsection{Matrices over semirings}
\label{ss:matrices}

Denote by $\Mat_{mn}(S)$ a set of all matrices $A = (a_{ij})$
with $m$~rows and $n$~columns whose coefficients belong to a semiring~$S$.
The sum $A \oplus B$ of matrices $A, B \in \Mat_{mn}(S)$ and the product
$AB$ of matrices $A \in \Mat_{lm}(S)$ and $B \in \Mat_{mn}(S)$ are defined
according to the usual rules of linear algebra:
$A\oplus B=(a_{ij} \oplus b_{ij})\in \mathrm{Mat}_{mn}(S)$ and
$$
AB=\left(\bigoplus_{k=1}^m a_{ij}\odot b_{kj}\right)\in\Mat_{ln}(S),
$$
where $A\in \Mat_{lm}(S)$ and $B\in\Mat_{mn}(S)$.
Note that we write $AB$ instead of $A\odot B$.

If the semiring~$S$ is
ordered, then the set $\Mat_{mn}(S)$ is ordered by the relation $A =
(a_{ij}) \preceq B = (b_{ij})$ iff $a_{ij} \preceq b_{ij}$ in~$S$ for all $1
\leqslant i \leqslant m$, $1 \leqslant j \leqslant n$.

The matrix multiplication is consistent with the order~$\preceq$ in the
following sense: if $A, A' \in \Mat_{lm}(S)$, $B, B' \in \Mat_{mn}(S)$ and
$A \preceq A'$, $B \preceq B'$, then $AB \preceq A'B'$ in $\Mat_{ln}(S)$. The set
$\Mat_{nn}(S)$ of square $(n \times n)$ matrices over an
idempotent semiring~$S$ forms a idempotent semiring with a
zero element $O = (o_{ij})$, where $o_{ij} = \0$, $1 \leqslant i, j
\leqslant n$, and a unit element $I = (\delta_{ij})$, where $\delta_{ij} =
\1$ if $i = j$ and $\delta_{ij} = \0$ otherwise.

The set $\Mat_{nn}$ is an example of a noncommutative semiring if $n>1$.

The closure operation in matrix semirings over an idempotent semiring~$S$ can
be defined inductively (another way to do that see in~\cite{Gol:99} and below): $A^*
= (a_{11})^* = (a^*_{11})$ in $\Mat_{11}(S)$ and for any integer $n > 1$
and any matrix
$$
   A = \begin{pmatrix} A_{11}& A_{12}\\ A_{21}& A_{22} \end{pmatrix},
$$
where $A_{11} \in \Mat_{kk}(S)$, $A_{12} \in \Mat_{k\, n - k}(S)$,
$A_{21} \in \Mat_{n - k\, k}(S)$, $A_{22} \in \Mat_{n - k\, n - k}(S)$,
$1 \leqslant k \leqslant n$, by defintion,
\begin{equation}
\label{A_Star}
   A^* = \begin{pmatrix}
   A^*_{11} \oplus A^*_{11} A_{12} D^* A_{21} A^*_{11} &
   \quad A^*_{11} A_{12} D^* \\[2ex]
   D^* A_{21} A^*_{11} &
   D^*
   \end{pmatrix},
\end{equation}
where $D = A_{22} \oplus A_{21} A^*_{11} A_{12}$. It can be proved that
this definition of $A^*$ implies that the equality $A^* = A^*A \oplus I$ is
satisfied and thus $A^*$ is a `regularized sum' of the series $I \oplus A
\oplus A^2 \oplus \dots$.

Note that this recurrence
relation coincides with the formulas of escalator method of matrix
inversion in the traditional linear algebra over the field of real or
complex numbers, up to the algebraic operations used. Hence this algorithm
of matrix closure requires a polynomial number of operations in~$n$.

\subsection{Discrete stationary Bellman equations}
\label{ss:bellman}

Let~$S$ be a semiring. The \emph{discrete stationary
Bellman equation} has the form
\begin{equation}
\label{AX+B}
	X = AX \oplus B,
\end{equation}
where $A \in \Mat_{nn}(S)$, $X, B \in \Mat_{ns}(S)$, and the matrix~$X$ is
unknown. Let $A^*$ be the closure of the matrix~$A$. It follows from the
identity $A^* = A^*A \oplus I$ that the matrix $A^*B$ satisfies this
equation; moreover, it can be proved that for idempotent semirings this
solution is the least in the set of solutions to equation~\eqref{AX+B}
with
respect to the partial order in $\Mat_{ns}(S)$.

Equation~\eqref{AX+B} over max-plus semiring arises in connection with
Bellman optimality principle and discretization of Hamilton-Jacobi equations, see e.g., \cite{McE:10}.
It is also intimately related with optimization problems
on graphs to be discussed below.

\subsection{Weighted directed graphs and matrices over semirings}

\begin{figure}[b]
\includegraphics[scale=0.8]{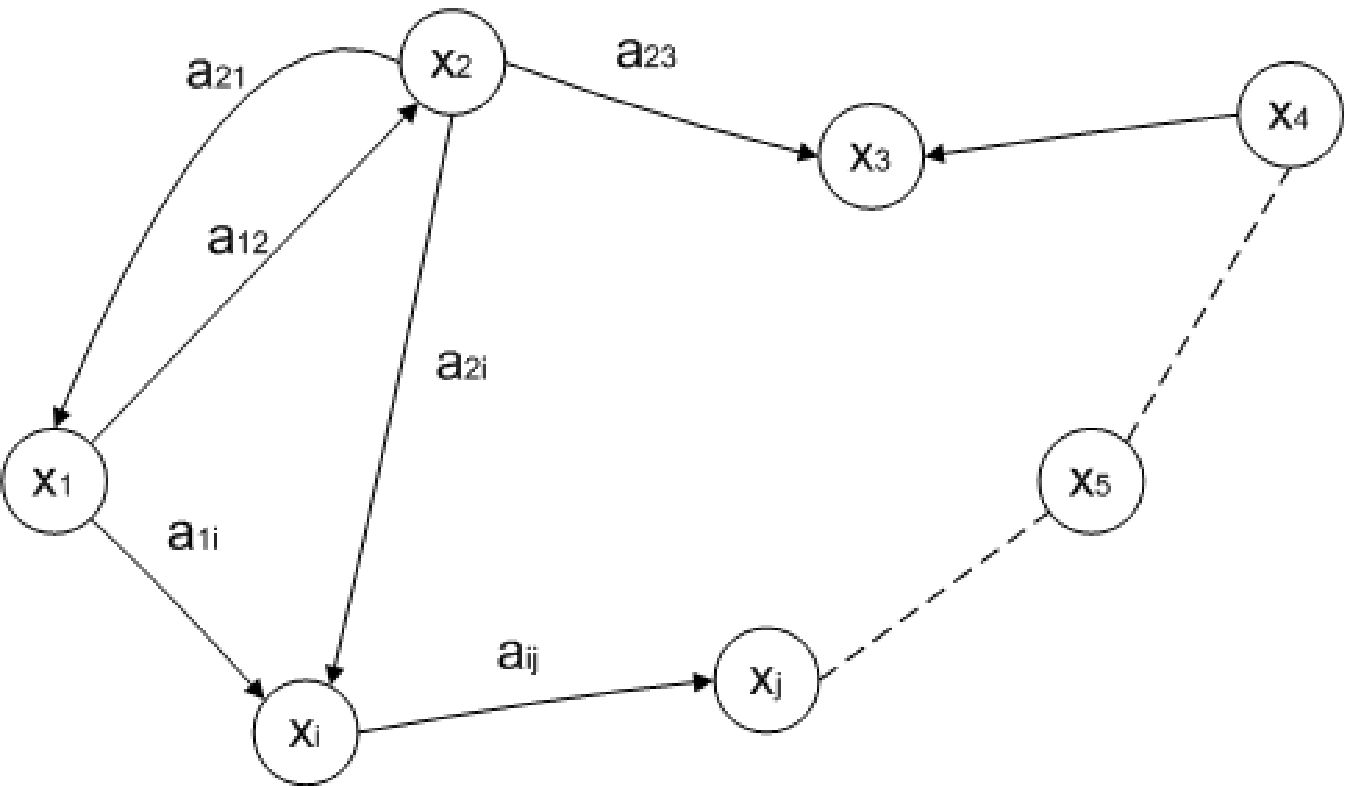}
%
%
\caption{A weighted directed graph.}
\label{f:graph}       
\end{figure}

Suppose that $S$ is a semiring with zero~$\0$ and unity~$\1$. It is well-known
that any square matrix $A = (a_{ij}) \in \Mat_{nn}(S)$ specifies a
\emph{weighted directed graph}. This geometrical construction includes
three kinds of objects: the set $X$ of $n$ elements $x_1, \dots, x_n$
called \emph{nodes}, the set $\Gamma$ of all ordered pairs $(x_i, x_j)$
such that $a_{ij} \neq \0$ called \emph{arcs}, and the mapping $A \colon
\Gamma \to S$ such that $A(x_i, x_j) = a_{ij}$. The elements $a_{ij}$ of
the semiring $S$ are called \emph{weights} of the arcs. See Fig.~\ref{f:graph}

Conversely, any given weighted directed graph with $n$ nodes specifies a
unique matrix $A \in \Mat_{nn}(S)$.

This definition allows for some pairs of nodes to be disconnected if the
corresponding element of the matrix $A$ is $\0$ and for some channels to be
``loops'' with coincident ends if the matrix $A$ has nonzero diagonal
elements. This concept is convenient for analysis of parallel and
distributed computations and design of computing media and networks (see,
e.g., \cite{Belov, MV:88, Voevod,LMRS}).

Recall that a sequence of nodes of the form
$$
	p = (y_0, y_1, \dots, y_k)
$$
with $k \geqslant 0$ and $(y_i, y_{i + 1}) \in \Gamma$, $i = 0, \dots, k -
1$, is called a \emph{path} of length $k$ connecting $y_0$ with $y_k$.
Denote the set of all such paths by $P_k(y_0,y_k)$. The weight $A(p)$ of a
path $p \in P_k(y_0,y_k)$ is defined to be the product of weights of arcs
connecting consecutive nodes of the path:
$$
	A(p) = A(y_0,y_1) \odot \cdots \odot A(y_{k - 1},y_k).
$$
By definition, for a `path' $p \in P_0(x_i,x_j)$ of length $k = 0$ the
weight is $\1$ if $i = j$ and $\0$ otherwise.

For each matrix $A \in \Mat_{nn}(S)$ define $A^0 = I = (\delta_{ij})$
(where $\delta_{ij} = \1$ if $i = j$ and $\delta_{ij} = \0$ otherwise) and
$A^k = AA^{k - 1}$, $k \geqslant 1$.  Let $a^{(k)}_{ij}$ be the $(i,j)$th
element of the matrix $A^k$. It is easily checked that
$$
   a^{(k)}_{ij} =
   \bigoplus_{\substack{i_0 = i,\, i_k = j\\
	1 \leqslant i_1, \ldots, i_{k - 1} \leqslant n}}
	a_{i_0i_1} \odot \dots \odot a_{i_{k - 1}i_k}.
$$
Thus $a^{(k)}_{ij}$ is the supremum of the set of weights corresponding to
all paths of length $k$ connecting the node $x_{i_0} = x_i$ with $x_{i_k} =
x_j$.

Denote the elements of the matrix $A^*$ by $a^{(*)}_{ij}$, $i, j = 1,
\dots, n$; then
$$
	a^{(*)}_{ij}
	= \bigoplus_{0 \leqslant k < \infty}
	\bigoplus_{p \in P_k(x_i, x_j)} A(p).
$$

The closure matrix $A^*$ solves the well-known \emph{algebraic path
problem}, which is formulated as follows: for each pair $(x_i,x_j)$
calculate the supremum of weights of all paths (of arbitrary length)
connecting node $x_i$ with node $x_j$. The closure operation in matrix
semirings has been studied extensively (see, e.g.,
\cite{AU:73,AHU:76, BCOQ, BJJM, BC-75, Car-71, Car:79, CG:79,CG:95,CG-91,
Gol:99,Gon-75,GM:79,GM:10,Gun:98, Kol-01, KM:97, LS-01}
and references therein).

\begin{example}[The shortest path problem.]
\label{ex:shortpath}
Let $S = \rmin$, so the weights are real numbers. In this case
$$
	A(p) = A(y_0,y_1) + A(y_1,y_2) + \dots + A(y_{k - 1},y_k).
$$
If the element $a_{ij}$ specifies the length of the arc $(x_i,x_j)$ in some
metric, then $a^{(*)}_{ij}$ is the length of the shortest path connecting
$x_i$ with $x_j$.
\end{example}

\begin{example}[The maximal path width problem.]
\label{ex:pathwidth}
Let $S = \rset \cup \{\0,\1\}$ with $\oplus = \max$, $\odot = \min$. Then
$$
	a^{(*)}_{ij} =
	\max_{p \in \bigcup\limits_{k \geqslant 1} P_k(x_i,x_j)} A(p),
	\quad
	A(p) = \min (A(y_0,y_1), \dots, A(y_{k - 1},y_k)).
$$
If the element $a_{ij}$ specifies the ``width'' of the arc
$(x_i,x_j)$, then the width of a path $p$ is defined as the minimal
width of its constituting arcs and the element $a^{(*)}_{ij}$ gives the
supremum of possible widths of all paths connecting $x_i$ with $x_j$.
\end{example}

\begin{example}[A simple dynamic programming problem.]
\label{ex:dinprog}
Let $S = \rmax$ and suppose $a_{ij}$ gives the \emph{profit} corresponding
to the transition from $x_i$ to $x_j$. Define the vector $B  = (b_i) \in
\Mat_{n1}(\rmax)$ whose element $b_i$ gives the \emph{terminal profit}
corresponding to exiting from the graph through the node $x_i$. Of course,
negative profits (or, rather, losses) are allowed. Let $m$ be the total
profit corresponding to a path $p \in P_k(x_i,x_j)$, i.e.
$$
	m = A(p) + b_j.
$$
Then it is easy to check that the supremum of profits that can be achieved
on paths of length $k$ beginning at the node $x_i$ is equal to $(A^kB)_i$
and the supremum of profits achievable without a restriction on the length
of a path equals $(A^*B)_i$.
\end{example}

\begin{example}[The matrix inversion problem.]
\label{ex:matrixinv}
Note that in the formulas of this section we are using distributivity of
the multiplication $\odot$ with respect to the addition $\oplus$ but do not
use the idempotency axiom. Thus the algebraic path problem can be posed for
a nonidempotent semiring $S$ as well (see, e.g.,~\cite{Rot-85}). For
instance, if $S = \rset$, then
$$
	A^* = I + A + A^2 + \dotsb = (I - A)^{-1}.
$$
If $\|A\| > 1$ but the matrix $I - A$ is invertible, then this expression
defines a regularized sum of the divergent matrix power series
$\sum_{i \geqslant 0} A^i$.
\end{example}

There are many other important examples of problems (in different areas) related to algorithms of linear algebra over semirings (transitive closures of relations, accessible sets, critical paths, paths of greatest capacities, the most reliable paths, interval and other problems), see
\cite{AU:73,AHU:76, Belov,BCOQ, BC-75, But:10, Car-71, Car:79, CC-02,
CGQ-99, CG:79,CG:95,CG-91, CGB-03, Fie+06,
Gol:99,Gon-75,GM:79,GM:10,Gun:98, Har+09,
Kol-01, KM:97, LS-00, LS-01, MV:88, Mys-05,Mys-06, Qua-90,MPlus-Scilab,RT-87,Rot-85,
Sed-92, Sim-88, Vor-63,Vor-67, Vor-70, Zim-06}.

We emphasize that this connection between the matrix closure operation and
solution to the Bellman equation gives rise to a number of different
algorithms for numerical calculation of the closure matrix. All these
algorithms are adaptations of the well-known algorithms of the traditional
computational linear algebra, such as the Gauss--Jordan elimination, various
iterative and escalator schemes, etc. This is a special case of the idempotent superposition principle.

In fact, the theory of the discrete stationary Bellman equation can be
developed using the identity $A^* = AA^* \oplus I$ as an additional axiom
without any substantial interpretation (the so-called \emph{closed
semirings}, see, e.g., \cite{BC-75, Gol:99, Leh-77, Rot-85}).

\section{Universal algorithms}
\label{s:universal}

Computational algorithms are constructed
on the basis of certain primitive operations. These operations manipulate
data that describe ``numbers.'' These ``numbers'' are elements of a
``numerical domain,'' i.e., a mathematical object such as the field of
real numbers, the ring of integers, or an idempotent semiring of numbers.

In practice elements of the numerical domains are replaced
by their computer representations, i.e., by elements of certain finite
models of these domains. Examples of models that can be conveniently used
for computer representation of real numbers are provided by various
modifications of floating point arithmetics, approximate arithmetics of
rational numbers~\cite{LRT-08}, and interval arithmetics. The difference
between mathematical objects (``ideal'' numbers) and their finite
models (computer representations) results in computational (e.g.,
rounding) errors.

An algorithm is called {\it universal\/} if it is independent of a
particular numerical domain and/or its computer representation.
A typical example of a universal algorithm is the computation of the
scalar product $(x,y)$ of two vectors $x=(x_1,\dots,x_n)$ and
$y=(y_1,\dots,y_n)$ by the formula $(x,y)=x_1y_1+\dots+x_ny_n$.
This algorithm (formula) is independent of a particular
domain and its computer implementation, since the formula is
well-defined for any semiring. It is clear that one algorithm can be
more universal than another. For example, the simplest Newton--Cotes formula, the
rectangular rule, provides the most universal algorithm for
numerical integration; indeed, this formula is valid even for
idempotent integration (over any idempotent semiring, see above and \cite{Belov, KM:97, Lit-07,LM-95,LM-98,LM:05,LMa-00,LMR-00, Mas-86, Mas-87a,Mas-87b,Mas:87}.
Other quadrature formulas (e.g., combined trapezoid rule or the Simpson
formula) are independent of computer arithmetics and can be
used (e.g., in an iterative form) for computations with
arbitrary accuracy. In contrast, algorithms based on
Gauss--Jacobi formulas are designed for fixed accuracy computations:
they include constants (coefficients and nodes of these formulas)
defined with fixed accuracy. Certainly, algorithms of this type can
be made more universal by including procedures for computing the
constants; however, this results in an unjustified complication of the
algorithms.

Computer algebra algorithms used in such systems as Mathematica,
Maple, REDUCE, and others are highly universal. Most of the standard
algorithms used in linear algebra can be rewritten in such a way
that they will be valid over any field and complete idempotent
semiring (including semirings of intervals; see below and~\cite{LS-00, LS-01, Sob-99}, where
an interval version of the idempotent linear algebra and the
corresponding universal algorithms are discussed).

As a rule, iterative algorithms (beginning with the successive approximation
method) for solving differential equations (e.g., methods of
Euler, Euler--Cauchy, Runge--Kutta, Adams, a number of important
versions of the difference approximation method, and the like),
methods for calculating elementary and some special functions based on
the expansion in Taylor's series and continuous fractions
(Pad\'e approximations) and others are independent of the computer
representation of numbers.

Calculations on computers usually are based on a floating-point arithmetic
with a mantissa of a fixed length; i.e., computations are performed
with fixed accuracy. Broadly speaking, with this approach only
the relative rounding error is fixed, which can lead to a drastic
loss of accuracy and invalid results (e.g., when summing series and
subtracting close numbers). On the other hand, this approach provides
rather high speed of computations. Many important numerical algorithms
are designed to use floating-point arithmetic (with fixed accuracy)
and ensure the maximum computation speed. However, these algorithms
are not universal. The above mentioned Gauss--Jacobi quadrature formulas,
computation of elementary and special functions on the basis of the
best polynomial or rational approximations or Pad\'e--Chebyshev
approximations, and some others belong to this type. Such algorithms
use nontrivial constants specified with fixed accuracy.

Recently, problems of accuracy, reliability, and authenticity of
computations (including the effect of rounding errors) have gained
much attention; in part, this fact is related to the ever-increasing
performance of computer hardware. When errors in initial data and
rounding errors strongly affect the computation results, such as in ill-posed
problems, analysis of stability of solutions, etc., it is often useful
to perform computations with improved and variable accuracy. In
particular, the rational arithmetic, in which the rounding error is
specified by the user~\cite{LRT-08}, can be used for this purpose.
This arithmetic is a useful complement to the interval analysis~\cite{Matij}.
The corresponding computational algorithms must be
universal (in the sense that they must be independent of the computer
representation of numbers).

\section{Universal algorithms of linear algebra over semirings}
\label{s:universal-linalg}
The most important linear algebra problem is to solve the system
of linear equations
\begin{equation}
\label{AX}
AX = B,
\end{equation}
where $A$ is a matrix with elements from the basic field and $X$ and
$B$ are vectors (or matrices) with elements from the same field.
It is required to find $X$ if $A$ and $B$ are given. If $A$ in~\eqref{AX}
is not the identity matrix $I$, then
system~\eqref{AX} can be written in form~\eqref{AX+B}, i.e.,
\begin{equation}
\label{AX+BB}
X = AX + B.
\end{equation}
It is well known that the form~\eqref{AX+BB} is convenient for using the
successive approximation method. Applying this method with the initial
approximation $X_0=0$, we obtain the solution
\begin{equation}
\label{A*B}
X = A^*B,
\end{equation}
where
\begin{equation}
\label{A*}
A^* = I+A+A^2+\cdots + A^n+\cdots
\end{equation}
On the other hand, it is clear that
\begin{equation}
\label{I-A-1}
A^* = (I-A)^{-1},
\end{equation}
if the matrix $I-A$ is invertible. The inverse matrix $(I-A)^{-1}$
can be considered as a regularized sum of the formal series~\eqref{A*}.

The above considerations can be extended to a broad class of
semirings.

The closure operation for matrix semirings ${\Mat}_n(S)$ can be defined
and computed in terms of the closure operation for $S$ (see Subsection~\ref{ss:matrices} above); some
methods are described in
\cite{AU:73,AHU:76,BC-75,Car-71,Car:79,Gol:99,Gon-75,GM:79,GM:10,KM:97,Kung-85,
LMa-00,LS-01, RT-87, Rot-85, Sed-92}.
 One such method is
described below ($LDM$-factorization), see \cite{LMRS}.

If $S$ is a field, then, by definition, $x^*=(1-x)^{-1}$ for any $x\ne 1$. If $S$ is an idempotent semiring, then, by definition,
\begin{equation}
x^*=\1\oplus x \oplus x^2 \oplus\cdots=\sup\{\1, x, x^2, \dots\},
\end{equation}
if this supremum exists. Recall that it exists if $S$ is complete, see section~4.2.

Consider a nontrivial universal algorithm applicable to matrices over
semirings with the closure operation defined.

\begin{example}[Semiring $LDM$-Factorization]
\label{ss:LDM}
Factorization of a matrix into the product $A = LDM$, where $L$ and $M$
are lower and upper triangular matrices with a unit diagonal,
respectively, and $D$ is a diagonal matrix, is used for solving
matrix equations $AX = B$. We construct a similar
decomposition for the Bellman equation $X = AX \oplus B$.

For the case $AX = B$, the decomposition $A = LDM$ induces the following
decomposition of the initial equation:
\begin{equation}
\label{LDM-dec}
   LZ = B, \qquad DY = Z, \qquad MX = Y.
\end{equation}
Hence, we have
\begin{equation}
   A^{-1} = M^{-1}D^{-1}L^{-1},
\label{AULinv}
\end{equation}
if $A$ is invertible. In essence, it is sufficient to find the
matrices $L$, $D$ and $M$, since the linear system~\eqref{LDM-dec} is easily
solved by a combination of the forward substitution for $Z$, the
trivial inversion of a diagonal matrix for $Y$, and the back
substitution for $X$.

Using~\eqref{LDM-dec} as a pattern, we can write
\begin{equation}
 \label{LDM}
  Z = LZ \oplus B, \qquad Y = DY \oplus Z, \qquad X = MX \oplus Y.
\end{equation}
Then
\begin{equation}
 \label{AMDLstar}
  A^* = M^*D^*L^*.
\end{equation}
A triple $(L,D,M)$ consisting of a lower triangular, diagonal, and
upper triangular matrices is called an $LDM$-{\it factorization} of a
matrix $A$ if relations~\eqref{LDM} and~\eqref{AMDLstar} are satisfied. We note that
in this case, the principal diagonals of $L$ and $M$ are zero.

The modification of the notion of $LDM$-factorization used in matrix
analysis for the equation $AX=B$ is constructed in analogy with a
construction suggested by Carr\'e in~\cite{Car-71,Car:79} for $LU$-factorization.

We stress that the algorithm described below can be applied to matrix
computations over any semiring under the condition that the unary
operation $a\mapsto a^*$ is applicable every time it is encountered
in the computational process. Indeed, when constructing the
algorithm, we use only the basic semiring operations of addition
$\oplus$ and multiplication $\odot$ and the properties of
associativity, commutativity of addition, and distributivity of
multiplication over addition.

If $A$ is a symmetric matrix over a semiring with a commutative
multiplication, the amount of computations can be halved, since
$M$ and $L$ are mapped into each other under transposition.

We begin with the case of a triangular matrix $A = L$ (or $A = M$).
Then, finding $X$ is reduced to the forward (or back) substitution.

\begin{center}
{\it Forward substitution}
\end{center}

 We are given:
\begin{itemize}
\item $L = \|l^i_j\|^n_{i,j = 1}$, where $l^i_j = \0$ for $i \leq j$
(a lower triangular matrix with a zero diagonal);
\item $B = \|b^i\|^n_{i = 1}$.
\end{itemize}

It is required to find the solution $X = \|x^i\|^n_{i = 1}$ to the
equation $X = LX \oplus B$. The program fragment solving this problem is as
follows.

\begin{tabbing}
   \qquad\=\qquad\=\kill
   for $i = 1$ to $n$ do\\*
   \{\> $x^i := b^i$;\\
   \> for $j = 1$ to $i - 1$ do\\*
   \>\> $x^i := x^i \oplus (l^i_j \odot x^j)$;\, \}\\
\end{tabbing}

\begin{center}
{\it Back substitution}
\end{center}

We are given
\begin{itemize}
\item $M = \|m^i_j\|^n_{i,j = 1}$, where $m^i_j = \0$ for $i \geq j$ (an
upper triangular matrix with a zero diagonal);
\item $B = \|b^i\|^n_{i = 1}$.
\end{itemize}

It is required to find the solution $X = \|x^i\|^n_{i = 1}$ to the
equation $X = MX \oplus B$. The program fragment solving this problem
is as follows.

\begin{tabbing}
   \qquad\=\qquad\=\kill
   for $i = n$ to 1 step $-1$ do\\*
   \{\> $x^i :=  b^i$;\\
   \> for $j = n$ to $i + 1$ step $-1$ do\\*
   \>\> $x^i :=  x^i \oplus (m^i_j \odot x^i)$;\, \}\\
\end{tabbing}

Both algorithms require $(n^2 - n) / 2$ operations $\oplus$ and $\odot$.

\begin{center}
{\it Closure of a diagonal matrix}
\end{center}

We are given
\begin{itemize}
\item $D = {\rm{diag}}(d_1, \ldots, d_n)$;
\item $B = \|b^i\|^n_{i = 1}$.
\end{itemize}

It is required to find the solution $X = \|x^i\|^n_{i = 1}$ to the
equation $X = DX \oplus B$. The program fragment solving this problem
is as follows.

\begin{tabbing}
   \qquad\=\qquad\=\kill
   for $i = 1$ to $n$ do\\*
   \> $x^i :=  (d_i)^* \odot b^i$;\\
\end{tabbing}

This algorithm requires $n$ operations $*$ and $n$ multiplications $\odot$.

\begin{center}
{\it General case}
\end{center}

We are given

\begin{itemize}
\item $L = \|l^i_j\|^n_{i,j = 1}$, where $l^i_j = \0$ if $i \leq j$;
\item $D = {\rm{diag}}(d_1, \ldots, d_n)$;
\item $M = \|m^i_j\|^n_{i,j = 1}$, where $m^i_j = \0$ if $i \geq j$;
\item $B = \|b^i\|^n_{i = 1}$.
\end{itemize}

It is required to find the solution $X = \|x^i\|^n_{i = 1}$ to the
equation $X = AX \oplus B$, where $L$, $D$, and $M$ form the
$LDM$-factorization of $A$. The program fragment solving this problem
is as follows.

\begin{tabbing}
        {\sc {FORWARD SUBSTITUTION}}\\*
   for $i = 1$ to $n$ do\\*
   \{\, $x^i :=  b^i$;\\*
   \, for $j = 1$ to $i - 1$ do\\*
   \,\, $x^i :=  x^i \oplus (l^i_j \odot x^j)$;\, \}\\
	\sc{CLOSURE OF A DIAGONAL MATRIX}\\*
   for $i = 1$ to $n$ do\\*
   \, $x^i :=  (d_i)^* \odot b^i$;\\
	\sc{BACK SUBSTITUTION}\\*
   for $i = n$ to 1 step $-1$ do\\*
   \{\, for $j = n$ to $i + 1$ step $-1$ do\\*
   \,\, $x^i :=  x^i \oplus (m^i_j \odot x^j)$;\, \}\\
\end{tabbing}

Note that $x^i$ is not initialized in the course of the back substitution.
The algorithm requires $n^2 - n$ operations $\oplus$, $n^2$ operations
$\odot$, and $n$ operations~$*$.

\begin{center}
{\it LDM-factorization}
\end{center}

We are given
\begin{itemize}
\item $A = \|a^i_j\|^n_{i,j = 1}$.
\end{itemize}

It is required to find the $LDM$-factorization of $A$:
$L = \|l^i_j\|^n_{i,j = 1}$, $D ={\rm{diag}}(d_1, \ldots, d_n)$, and
$M = \|m^i_j\|^n_{i,j = 1}$, where $l^i_j = \0$ if $i \leq j$, and
$m^i_j = \0$ if $i \geq j$.

The program uses the following internal variables:
\begin{itemize}
\item $C = \|c^i_j\|^n_{i,j = 1}$;
\item $V = \|v^i\|^n_{i = 1}$;
\item $d$.
\end{itemize}

\begin{tabbing}
   \qquad\=\qquad\=\qquad\=\kill
   \sc{INITIALISATION}\\*
	for $i = 1$ to $n$ do\\*
	\> for $j = 1$ to $n$ do\\*
	\>\> $c^i_j = a^i_j$;\\
	\sc{MAIN LOOP}\\*
	for $j = 1$ to $n$ do\\*
	\{\> for $i = 1$ to $j$ do\\*
	\>\> $v^i :=  a^i_j$;\\
	\> for $k = 1$ to $j - 1$ do\\*
	\>\> for $i = k + 1$ to $j$ do\\*
	\>\>\> $v^i :=  v^i \oplus (a^i_k \odot v^k)$;\\
	\> for $i = 1$ to $j - 1$ do\\*
	\>\> $a^i_j :=  (a^i_i)^* \odot v^i$;\\
	\> $a^j_j :=  v^j$;\\
	\> for $k = 1$ to $j - 1$ do\\*
	\>\> for $i = j + 1$ to $n$ do\\*
	\>\>\> $a^i_j :=  a^i_j \oplus (a^i_k \odot v^k)$;\\
	\> $d = (v^j)^*$;\\
	\> for $i = j + 1$ to $n$ do\\*
	\>\> $a^i_j :=  a^i_j \odot d$;\, \}\\
\end{tabbing}

This algorithm requires $(2n^3 - 3n^2 + n) /6$ operations $\oplus$, $(2n^3 +
3n^2 -5n) / 6$ operations $\odot$, and $n(n + 1) / 2$ operations $*$.
After its completion, the matrices $L$, $D$, and $M$ are contained,
respectively, in the lower triangle, on the diagonal, and in the upper
triangle of the matrix $C$. In the case when $A$ is symmetric about the
principal diagonal and the semiring over which the matrix is defined
is commutative, the algorithm can be modified in such a way that the
number of operations is reduced approximately by a factor of two.
\end{example}

Other examples can be found in \cite{
Car-71,Car:79,Gol:99,Gon-75,GM:79,GM:10,Kung-85,Leh-77,Rot-85,Sed-92}.

Note that to compute the matrices $A^*$ and $A^*B$ it is convenient to solve the Bellman equation~(\ref{AX+BB}).

Some other interesting and important problems of linear algebra over
semirings are examined, e.g., in
\cite{BN-74,But:10,BZ-06,CC-02, CGB-03,Fie+06,Gol:99,Gon-75,GM:79,GM:10,Har+09,Mys-05,Mys-06,Neu:90,Pan-61,
Vor-63,Vor-67,Vor-70,Zim-06}.

\begin{remark}
It is well known that linear problems and equations are especially convenient for parallelization, see, e.g.,~\cite{Voevod}.
Standard methods (including the so-called block methods) constructed in the framework of the traditional mathematics can be extended to universal algorithms over semirings (the correspondence principle!). For example, formula~\eqref{A_Star} discussed in Subsection~\ref{ss:matrices} leads to a simple block method for parallelization of the closure operations. Other standard methods of linear algebra~\cite{Voevod} can be used in a similar way.
\end{remark}

\section{The correspondence principle for computations}
\label{s:corresp-comp}
Of course, the idempotent correspondence principle is valid for
algorithms as well as for their software and hardware
implementations~\cite{LM-95,LM-98,LMa-00,LMR-00}. Thus:

{\it If we have an important and interesting numerical algorithm, then
there is a good chance that its semiring analogs are important and
interesting as well.}

In particular, according to the superposition principle,
analogs of linear
algebra algorithms are especially important. Note that
numerical algorithms
for standard infinite-dimensional linear problems over idempotent
semirings (i.e., for
problems related to idempotent integration, integral operators and
transformations, the Hamilton-Jacobi and generalized Bellman equations)
deal with the corresponding finite-dimensional (or finite) ``linear
approximations''. Nonlinear algorithms often can be approximated by linear
ones. Thus the idempotent linear algebra is a basis for the idempotent
numerical analysis.

Moreover, it is well-known that linear algebra algorithms easily lend themselves to parallel computation; their idempotent analogs admit
parallelization as
well. Thus we obtain a systematic way of applying parallel computing to
optimization problems.

Basic algorithms of linear algebra (such as inner product of two vectors,
matrix addition and multiplication, etc.) often do not depend on
concrete semirings, as well as on the nature of domains containing the
elements of vectors and matrices. Algorithms to construct the closure
$A^*=I\oplus A\oplus A^2\oplus\cdots\oplus A^n\oplus\cdots=
\bigoplus^{\infty}_{n=1} A^n$ of an idempotent matrix $A$ can be derived
from standard methods for calculating $(I-A)^{-1}$. For the
Gauss--Jordan elimination method (via LU-decomposition) this trick was used in~\cite{Rot-85},
and the corresponding algorithm is universal and can be applied both to
the Bellman equation and to computing the inverse of a real (or complex)
matrix $(I - A)$. Computation of $A^{-1}$ can be derived
from this universal
algorithm with some obvious cosmetic transformations.

Thus it seems reasonable to develop universal algorithms that can deal
equally well with initial data of different domains sharing the same
basic structure~\cite{LM-95,LM-98,LMR-00}.

\section{The correspondence principle for hardware design}
\label{s:corresp-hard}
A systematic application of the correspondence principle to computer
calculations leads to a unifying approach to software and hardware
design.

The most important and standard numerical algorithms have many hardware
realizations in the form of technical devices or special processors.
{\it These devices often can be used as prototypes for new hardware
units generated by substitution of the usual arithmetic operations
for its semiring analogs and by addition tools for performing neutral
elements $\0$ and} $\1$ (the latter usually is not difficult). Of course,
the case of numerical semirings consisting of real numbers (maybe except
neutral elements)  and semirings of numerical intervals is the most simple and natural
\cite{Lit-07,LM-95,LM-98,LM:05,LMa-00,LMR-00,LS-00,LS-01,Sob-99}.
Note that for semifields (including $\rmax$ and $\rmin$)
the operation of division is also defined.

Good and efficient technical ideas and decisions can be transferred
from prototypes to new hardware units. Thus the correspondence
principle generated a regular heuristic method for hardware design.
Note that to get a patent it is necessary to present the so-called
`invention formula', that is to indicate a prototype for the suggested
device and the difference between these devices.

Consider (as a typical example) the most popular and important algorithm
of computing the scalar product of two vectors:
\begin{equation}
\label{scalprod}
(x,y)=x_1y_1+x_2y_2+\cdots + x_ny_n.
\end{equation}
The universal version of~\eqref{scalprod} for any semiring $A$ is obvious:
\begin{equation}
\label{univscal}
(x,y)=(x_1\odot y_1)\oplus(x_2\odot y_2)\oplus\cdots\oplus
(x_n\odot y_n).
\end{equation}
In the case $A=\rmax$ this formula turns into the following one:
\begin{equation}
\label{maxscal}
(x,y)=\max\{ x_1+y_1,x_2+y_2, \cdots, x_n+y_n\}.
\end{equation}

This calculation is standard for many optimization algorithms, so
it is useful to construct a hardware unit for computing~\eqref{maxscal}. There
are many different devices (and patents) for computing~\eqref{scalprod} and every
such device can be used as a prototype to construct a new device for
computing~\eqref{maxscal} and even~\eqref{univscal}. Many processors for matrix multiplication
and for other algorithms of linear algebra are based on computing
scalar products and on the corresponding ``elementary'' devices
respectively, etc.

There are some methods to make these new devices more universal than
their prototypes. There is a modest collection of possible operations
for standard numerical semirings: max, min, and the usual arithmetic
operations. So, it is easy to construct programmable hardware
processors with variable basic operations. Using modern technologies
it is possible to construct cheap special-purpose multi-processor
chips implementing examined algorithms. The so-called
systolic processors are
especially convenient for this purpose. A systolic array is a
`homogeneous' computing medium consisting of elementary
processors, where the general scheme and processor connections
are simple and regular. Every elementary processor pumps data in and
out performing elementary operations in a such way that the
corresponding data flow is kept up in the computing medium; there
is an analogy with the blood circulation and this is a reason for the
term ``systolic'', see e.g.,
\cite{LM-95,LM-98, LMR-00, LMRS, Mas-tech91,RT-87,Rot-85,Sed-92}.

Some systolic processors for the general algebraic path problem are
presented in~\cite{RT-87,Rot-85,Sed-92}. In particular, there is a systolic array of
$n(n+1)$ elementary processors which performs computations of the Gauss--Jordan
elimination algorithm and can solve the algebraic path problem within $5n-2$
time steps. Of course, hardware implementations for important and popular basic
algorithms increase the speed of data processing.

The so-called GPGPU (General-Purpose computing on Graphics Processing Units) technique is another important field for applications of the correspondence principle. The matter is that graphic processing units (hidden in modern laptop and desktop computers) are potentially powerful processors for solving numerical problems. The recent tremendous progress in graphical processing hardware and software resulted in new ``open'' programmable parallel computational devices (special processors), see, e.g.,~\cite{IEEE-09,Bli-08,Owe-08}. These devices are going to be standard for coming PC (personal computers) generations. Initially used for graphical processing only (at that time they were called GPU), today they are used for various fields, including audio and video processing, computer simulation, and encryption.   But this list can be considerably enlarged following the correspondence principle: the basic operations would be used as parameters. Using the technique described in this paper (see also our references), standard linear algebra algorithms can be used for solving different problems in different areas. In fact, the hardware supports all operations needed for the most important idempotent semirings: plus, times, min, max. The most popular linear algebra packages [ATLAS (Automatically Tuned Linear Algebra Software), LAPACK, PLASMA (Parallel Linear Algebra for Scalable Multicore Architectures)] can already use GPGPU, see \cite{ATLAS,LAPACK,PLASMA}. We propose to make these tools more powerful by using parameterized algorithms.

Linear algebra over the most important numerical semirings generates solutions for many concrete problems in different areas, see above.

Note that to be consistent with operations we have to redefine zero (0) and unit (1) elements (see above); comparison operations must be also redefined as it is described above. Once the operations are redefined, then the most of basic linear algebra algorithms, including back and forward substitution, Gauss elimination method, Jordan elimination method and others could be rewritten for new domains and data structures.
Combined with the power of the new parallel hardware this approach could change PC from entertainment devices to power full instruments.

\section{The correspondence principle for software design}
\label{s:corresp-soft}

Software implementations for universal semiring algorithms are not
as efficient as hardware ones (with respect to the computation speed)
but they are much more flexible. Program modules can deal with abstract (and
variable) operations and data types. These
operations and data types can be defined by the corresponding
input data. In this case they can be generated
by means of additional program modules. For programs written in
this manner it is convenient to use special techniques of the
so-called object oriented (and functional) design, see, e.g.,
\cite{Lor:93, Pohl:97,SL:94}. Fortunately, powerful tools supporting the
object-oriented software design have recently appeared including compilers
for real and convenient programming languages (e.g. $C^{++}$ and Java) and modern computer algebra systems.

Recently, this type of programming technique has been dubbed
generic programming (see, e.g.,~\cite{BJJM, Pohl:97}). To help automate the
generic programming, the so-called Standard Template Library (STL)
was developed in the framework of $C^{++}$~\cite{Pohl:97,SL:94}. However,
high-level tools, such as STL, possess both obvious advantages
and some disadvantages and must be used with caution.

It seems that it is natural to obtain an implementation of the correspondence
principle approach to scientific calculations in the form of a
powerful software system based on a collection of universal
algorithms. This approach ensures a working time reduction for
programmers and users because of the software unification.
The arbitrary necessary accuracy and safety of numeric calculations can be ensured
as well.

This software system may be especially useful for designers
of algorithms, software engineers, students and mathematicians.

Note that there are some software systems oriented to calculations with idempotent semirings like $\rmax$; see, e.g.,~\cite{MPlus-Scilab}. However these systems do not support universal algorithms.

\section{Interval analysis in idempotent mathematics}
\label{s:interval}

Traditional interval analysis is a nontrivial and popular mathematical area, see, e.g.,~\cite{AH:83,Fie+06,Kre+98,Matij,Moo:79,Neu:90}. An ``idempotent'' version of interval analysis (and moreover interval analysis over positive semirings) appeared in~\cite{LS-00,LS-01,Sob-99}. Later the idempotent interval analysis has attracted
many experts in tropical linear algebra and applications, see, e.g.,~\cite{CC-02,Fie+06,Har+09,Mys-05,Mys-06, Zim-06}. We also mention
the closely related interval analysis over the positive semiring $\textbf{R}_+$ discussed in~\cite{BN-74}.

Let a set~$S$ be partially ordered by a relation $\preceq$.
A \emph{closed interval} in~$S$ is a subset of the form $\x = [\lx, \ux] =
\{\, x \in S \mid \lx \preceq x \preceq \ux\, \}$, where the elements $\lx \preceq
\ux$ are called \emph{lower} and \emph{upper bounds} of the interval $\x$.
The order~$\preceq$ induces a partial ordering on the set of all closed
intervals in~$S$: $\x \preceq \y$ iff $\lx \preceq \ly$ and $\ux \preceq \uy$.

A \emph{weak interval extension} $I(S)$ of an ordered semiring~$S$ is the
set of all closed intervals in~$S$ endowed with operations $\oplus$
and~$\odot$ defined as ${\x \oplus \y} = [{\lx \oplus \ly}, {\ux \oplus
\uy}]$, ${\x \odot \y} = [{\lx \odot \ly}, {\ux \odot \uy}]$ and a partial
order induced by the order in $S$. The closure operation in $I(S)$ is
defined by $\x^* = [\lx^*, \ux^*]$. There are some other interval extensions (including the so-called strong interval extension~\cite{LS-01}) but the weak extension is more convenient.

The extension $I(S)$ is idempotent if $S$ is an idempotent semiring.
A universal algorithm over $S$ can be applied to $I(S)$ and we shall get an interval version of the initial algorithm.
Usually both the versions have the same complexity. For the discrete stationary Bellman equation and the corresponding optimization problems on graphs, interval analysis was examined in~\cite{LS-00,LS-01} in details. Other problems of idempotent linear algebra were examined in~\cite{CC-02,Fie+06,Har+09,Mys-05,Mys-06}.

Idempotent mathematics appears to be remarkably simpler than its
traditional analog. For example, in traditional interval arithmetic,
multiplication of intervals is not distributive with respect to addition of
intervals, whereas in idempotent interval arithmetic this distributivity is
preserved. Moreover, in traditional interval analysis the set of all
square interval matrices of a given order does not form even a semigroup
with respect to matrix multiplication: this operation is not associative
since distributivity is lost in the traditional interval arithmetic. On the
contrary, in the idempotent (and positive) case associativity is preserved. Finally, in
traditional interval analysis some problems of linear algebra, such as
solution of a linear system of interval equations, can be very difficult
(more precisely, they are $NP$-hard, see~
\cite{Cox-99,Fie+06,Kre+93,Kre+98} and references therein). It was noticed  in~\cite{LS-00,LS-01} that in the idempotent case solving an interval linear system
requires a polynomial number of operations (similarly to the usual Gauss
elimination algorithm).  The remarkable simplicity of idempotent interval
arithmetic is due to the following properties: the monotonicity of arithmetic operations and the positivity of all elements of an idempotent semiring.

Interval estimates in idempotent mathematics are usually exact. In the traditional theory such estimates tend to be overly pessimistic.

\section*{Acknowledgement}
The author is sincerely grateful to V.~N.~Kolokoltsov, V.~P.~Maslov, S.~N.~Sergeev,
A.~N.~Sobolevski and A.~V.~Tchourkin for valuable suggestions, help and support.

This work is supported by the RFBR grants.


\begin{thebibliography}{100}

\bibitem{AHU:76}
A.V. Aho, J.E. Hopcroft, and J.D. Ullman.
\newblock {\em The Design and Analysis of Computer Algorithms}.
\newblock Addison Wesley Publ. Co., Reading, MS, 1976.

\bibitem{AU:73}
A.V. Aho and J.D. Ullman.
\newblock {\em The Theory of Parsing, Translation and Compiling. Vol. 2:
  Compiling}.
\newblock Prentice Hall, Englewood Cliffs, NJ, 1973.

\bibitem{AGK}
M.~Akian, S.~Gaubert, and V.~Kolokoltsov.
\newblock Set coverings and invertibility of the functional Galois connections.
\newblock In G.~Litvinov and V.~Maslov, editors, {\em Idempotent Mathematics
  and Mathematical Physics}, volume 377, pages 19--51. American Mathematical
  Society, Providence, 2005.
\newblock E-print arXiv:math.FA/0403441.

\bibitem{AH:83}
G.~Alefeld and J.~Herzberger.
\newblock {\em Introduction to Interval Computations}.
\newblock Academic Press, New York, 1983.

\bibitem{Belov}
S.M. Avdoshin, V.V. Belov, V.P. Maslov, and A.M. Chebotarev.
\newblock Design of computational media: mathematical aspects.
\newblock In V.P. Maslov and K.A. Volosov, editors, {\em Mathematical aspects
  of computer engineering}, pages 9--145. Mir Publishers, Moscow, 1988.

\bibitem{BCOQ}
F.~L. Baccelli, G.~Cohen, G.~J. Olsder, and J.~P. Quadrat.
\newblock {\em Synchronization and Linearity: an Algebra for Discrete Event
  Systems}.
\newblock Wiley, 1992.

\bibitem{BC-75}
R.~C. Backhouse and B.~A. Carr{\'{e}}.
\newblock Regular algebra applied to path-finding problems.
\newblock {\em J. of Inst. of Maths. and Applics}, 15:161--186, 1975.

\bibitem{BJJM}
R.C. Backhouse, P.~Janssen, J.~Jeuring, and L.~Meertens.
\newblock Generic programming - an introduction.
\newblock In {\em Lecture Notes in Comp. Sci.}, volume 1608, pages 28--115.
  1999.

\bibitem{BN-74}
W.~Barth and E.~Nuding.
\newblock Optimale L{\"o}sung von Intervalgleichungsystemen.
\newblock {\em Computing}, 12:117--125, 1974.

\bibitem{Bir:67}
G.~Birkhoff.
\newblock {\em Lattice Theory}.
\newblock Amer.Math.Soc., Providence, 1967.

\bibitem{Bli-08}
D.~Blithe.
\newblock Rise of the graphics processors.
\newblock {\em Proc. of the IEEE}, 96(5):761--778, 2008.

\bibitem{But:10}
P.~Butkovi{\v{c}}.
\newblock {\em Max-linear Systems: Theory and Algorithms}.
\newblock Springer, London, 2010.

\bibitem{BZ-06}
P.~Butkovi{\v{c}} and K.~Zimmermann.
\newblock A strongly polynomial algorithm for solving two-sided linear systems
  in max-algebra.
\newblock {\em Discrete Appl. Math.}, 154:437--446, 2006.

\bibitem{Car-71}
B.A. Carr\'{e}.
\newblock An algebra for network routing problems.
\newblock {\em J. of the Inst. of Maths. and Applics}, 7:273--294, 1971.

\bibitem{Car:79}
B.A. Carr{\'e}.
\newblock {\em Graphs and Networks}.
\newblock The Clarendon Press/Oxford Univ. Press, Oxford, 1979.

\bibitem{CC-02}
K.~Cechl{\'a}rov{\'a} and R.A. Cuninghame-Green.
\newblock Interval systems of max-separable linear equations.
\newblock {\em Linear Alg. Appl.}, 340(1-3):215--224, 2002.

\bibitem{CGQ-99}
G.~Cohen, S.~Gaubert, and J.~P. Quadrat.
\newblock Max-plus algebra and system theory: where we are and where to go now.
\newblock {\em Annual Reviews in Control}, 23:207--219, 1999.

\bibitem{CGQ-04}
G.~Cohen, S.~Gaubert, and J.~P. Quadrat.
\newblock Duality and separation theorems in idempotent semimodules.
\newblock {\em Linear Alg. Appl.}, 379:395--422, 2004.
\newblock E-print \arxiv{math.FA/0212294}.

\bibitem{Cox-99}
G.E. Coxson.
\newblock Computing exact bounds on the elements of an inverse interval matrix
  is NP-hard.
\newblock {\em Reliable Computing}, 5:137--142, 1999.

\bibitem{CG:79}
R.~A. Cuninghame-Green.
\newblock {\em Minimax Algebra}, volume 166 of {\em Lecture Notes in Economics
  and Mathematical Systems}.
\newblock Springer, Berlin, 1979.

\bibitem{CG-91}
R.~A. Cuninghame-Green.
\newblock Minimax algebra and its applications.
\newblock {\em Fuzzy Sets and Systems}, 41:251--267, 1991.

\bibitem{CG:95}
R.~A. Cuninghame-Green.
\newblock Minimax algebra and applications.
\newblock {\em Advances in Imaging and Electron Physics}, 90:1--121, 1995.

\bibitem{CGB-03}
R.~A. Cuninghame-Green and P.~Butkovi{\v{c}}.
\newblock The equation ${A}\otimes x={B}\otimes y$ over (max,+).
\newblock {\em Theoretical Computer Science}, 293:3--12, 2003.



\bibitem{Fie+06}
M.~Fiedler, J.~Nedoma, J.~Ram{\'{\i}}k, J.~Rohn, and K.~Zimmermann.
\newblock {\em Linear Optimization Problems with Inexact Data}.
\newblock Springer, New York, 2006.

\bibitem{GKZ}
I.~M. Gel�fand, M.~Kapranov, and A.~Zelevinsky.
\newblock {Multidimensional Determinants, Discriminants and Resultants}.
\newblock {\em Birkh{\"a}user, Boston}, 1994.

\bibitem{Gol:99}
J.~Golan.
\newblock {\em Semirings and Their Applications}.
\newblock Kluwer, 2000.

\bibitem{Gon-75}
M.~Gondran.
\newblock Path algebra and algorithms.
\newblock In B.~Roy, editor, {\em Combinatorial programming: methods and
  applications}, pages 137--148. Reidel, Dordrecht, 1975.

\bibitem{GM:79}
M.~Gondran and M.~Minoux.
\newblock {\em Graphes et Algorithmes}.
\newblock Editions Eylrolles, Paris, 1979.

\bibitem{GM:10}
M.~Gondran and M.~Minoux.
\newblock {\em Graphs, {D}ioids and {S}emirings. New {M}odels and
  {A}lgorithms}.
\newblock Springer, 2010.

\bibitem{Gun:98}
J.~Gunawardena, editor.
\newblock {\em {Idempotency}}. Volume 11 of Publ. of the I. Newton Institute, Cambridge Univ. Press, Cambridge, 1998.

\bibitem{Har+09}
L.~Hardouin, B.~Cottenceau, M.~Lhommeau, and E.~Le Corronc.
\newblock Interval systems over idempotent semiring.
\newblock {\em Linear Alg. Appl.}, 431:855--862, 2009.

\bibitem{IMS:07}
I.~Itenberg, G.~Mikhalkin, and E.~Shustin.
\newblock {\em Tropical {A}lgebraic {G}eometry}, volume~35 of {\em Oberwolfach
  {S}eminars}.
\newblock Birkh{\"a}user, Basel et al., 2007.

\bibitem{KM:97}
V.~N. Kolokoltsov and V.~P. Maslov.
\newblock {\em Idempotent Analysis and its Applications}.
\newblock Kluwer Academic Publ., 1997.

\bibitem{Kol-01}
V.N. Kolokoltsov.
\newblock Idempotency structures in optimization.
\newblock {\em J. of Math. Sci.}, 104(1):847--880, 2001.

\bibitem{Kre+98}
V.~Kreinovich, A.~Lakeev, J.~Rohn, and P.~Kahl.
\newblock {\em Computational Complexity and Feasibility of Data Processing and
  Interval Computations}.
\newblock Kluwer Academic Publishers, Dordrecht, 1998.

\bibitem{Kre+93}
V.~Kreinovich, A.V. Lakeyev, and S.I. Noskov.
\newblock Optimal solution of interval systems is intractable (NP-hard).
\newblock {\em Interval computations}, 1:6--14, 1993.

\bibitem{Kung-85}
H.T. Kung.
\newblock Two-level pipelined systolic arrays for matrix multiplication,
  polynomial evaluation and discrete Fourier transformation.
\newblock In J.~Demongeof et~al., editor, {\em Dynamic and Cellular Automata},
  pages 321--330. Academic Press, New York et al., 1985.

\bibitem{Leh-77}
D.J. Lehmann.
\newblock Algebraic structures for transitive closure.
\newblock {\em Theoret. Comp. Sci.}, 4:59--76, 1977.









\bibitem{Lit-07}
G.~L. Litvinov.
\newblock The Maslov dequantization, idempotent and tropical mathematics: a
  brief introduction.
\newblock {\em Journal of Mathematical Sciences}, 140(3):426--441, 2007.
\newblock E-print \arxiv{math.GM/0507014}.

\bibitem{LM-95}
G.~L.~Litvinov and V.~P.~Maslov.
\newblock {The correspondence principle for idempotent calculus and some
  computer applications}.
\newblock In J.~Gunawardena, editor, {\em Idempotency}, pages 420--443.
  Cambridge Univ. Press, 1998.
\newblock E-print \arxiv{math/0101021}.

\bibitem{LM-96}
G.~L. Litvinov and V.~P. Maslov.
\newblock Idempotent mathematics: correspondence principle and applications.
\newblock {\em Russian Mathematical Surveys}, 51(6):1210--1211, 1996.

\bibitem{LM-98}
G.~L. Litvinov and V.~P. Maslov.
\newblock {\em Correspondence principle for idempotent calculus and some
  computer applications.}
\newblock (IHES/M/95/33), Institut des Hautes Etudes Scientifiques, Bures-sur-Yvette, 1995.
\newblock E-print \arxiv{math.GM/0101021}.



\bibitem{LM:05}
G.~L.~Litvinov and V.~P.~Maslov, editors.
\newblock {\em {Idempotent Mathematics and Mathematical Physics}}, volume 307
  of {\em Contemporary Mathematics}. Amer. Math. Soc., Providence, 2005.

\bibitem{LMR-00}
G.~L. Litvinov, V.P. Maslov, and A.Ya. Rodionov.
\newblock {\em A Unifying Approach to Software and Hardware Design for
  Scientific Calculations and Idempotent Mathematics}.
\newblock International Sophus Lie Centre, Moscow, 2000.
\newblock E-print \arxiv{math.SC/0101069}.

\bibitem{LMRS}
G.~L. Litvinov, V.P. Maslov, A.Ya. Rodionov, and A.N. Sobolevski.
\newblock Universal algorithms, mathematics of semirings and parallel
  computations.
\newblock {\em Lecture Notes in Computational Science and Engineering},
  75:63--89, 2011.
\newblock E-print \arxiv{1005.1252}.

\bibitem{LMS:07}
G. Litvinov, V. Maslov, and S. Sergeev, editors.
\newblock {\em {Idempotent and Tropical Mathematics and Problems of
  Mathematical Physics (Volumes I and II)}}, Moscow, 2007. French-Russian
  Laboratory J.V.~Poncelet.
\newblock E-prints \arxiv{0710.0377} and \arxiv{0709.4119}.


\bibitem{LMS-98}
G.~L. Litvinov, V.P. Maslov, and G.B. Shpiz.
\newblock Linear functionals on idempotent spaces: an algebraic approach.
\newblock {\em Doklady Mathematics}, 58(3):389--391, 1998.
\newblock E-print \arxiv{math.FA/0012268}.

\bibitem{LMS-99}
G.~L. Litvinov, V.P. Maslov, and G.B. Shpiz.
\newblock Tensor products of idempotent semimodules. {A}n algebraic approach.
\newblock {\em Math. Notes}, 65(4):497--489, 1999.
\newblock E-print \arxiv{math.FA/0101153}.

\bibitem{LMS-01}
G.~L. Litvinov, V.P. Maslov, and G.B. Shpiz.
\newblock Idempotent functional analysis: An algebraic approach.
\newblock {\em Math. Notes}, 69(5):758--797, 2001.

\bibitem{LMS-02}
G.~L. Litvinov, V.P. Maslov, and G.B. Shpiz.
\newblock Idempotent (asymptotic) analysis and the representation theory.
\newblock In V.A. Malyshev and A.M. Vershik, editors, {\em Asymptotic
  Combinatorics with Applications to Mathematical Physics}, pages 267--268.
  Kluwer Academic Publ., Dordrecht et al., 2002.

\bibitem{LMa-00}
G.~L. Litvinov and E.V. Maslova.
\newblock Universal numerical algorithms and their software implementation.
\newblock {\em Programming and Computer Software}, 26(5):275--280, 2000.
\newblock E-print \arxiv{math.SC/0102114}.

\bibitem{LRT-08}
G.~L. Litvinov, A.Ya. Rodionov, and A.V. Tchourkin.
\newblock Approximate rational arithmetics and arbitrary precision computations
  for universal algorithms.
\newblock {\em Internat. J. of Pure and Applied Math.}, 45(2):193--204, 2008.
\newblock E-print \arxiv{math.NA/0101152}.

\bibitem{LS:09}
G.~L. Litvinov and S.N. Sergeev, editors.
\newblock {\em {Tropical and Idempotent Mathenatics}}, volume 495 of {\em
  Contemporary Mathematics}. Amer.Math. Soc., Providence, 2009.

\bibitem{LSz-02}
G.~L. Litvinov and G.B. Shpiz.
\newblock Nuclear semimodules and kernel theorems in idempotent analysis: an
  algebraic approach.
\newblock {\em Doklady Mathematics}, 66(2):197--199, 2002.
\newblock E-print \arxiv{math.FA/0202386}.

\bibitem{LSz-05}
G.~L. Litvinov and G.B. Shpiz.
\newblock The dequantization transform and generalized Newton polytopes.
\newblock In G.~Litvinov and V.~Maslov, editors, {\em Idempotent Mathematics
  and Mathematical Physics}, volume 377 of {\em Contemporary Mathematics}, pages 181--186. American Mathematical
  Society, Providence, 2005.



\bibitem{LSz-07a}
G.~L. Litvinov and G.B. Shpiz.
\newblock The dequantization procedures related to maslov dequantization.
\newblock In G.~Litvinov, V.~Maslov, and S.~Sergeev, editors, {\em Idempotent
  and tropical mathematics and problems of mathematical physics (Volume I)},
  pages 99--104, 2007.
\newblock E-print \arxiv{0710.0377}.

\bibitem{LSz-07b}
G.~L. Litvinov and G.B. Shpiz.
\newblock Kernel theorems and nuclearity in idempotent mathematics. An
  algebraic approach.
\newblock {\em Journal of Mathematical Sciences.}, 141(4):1417--1428, 2007.
\newblock E-print \arxiv{math.FA/0609033}.




\bibitem{LS-00}
G.~L. Litvinov and A.N. Sobolevski{\u{\i}}.
\newblock Exact interval solutions of the discrete bellman equation and
  polynomial complexity in interval idempotent linear algebra.
\newblock {\em Doklady Mathematics}, 62(2):199--201, 2000.
\newblock E-print \arxiv{math.LA/0101041}.

\bibitem{LS-01}
G.~L. Litvinov and A.N. Sobolevski{\u{\i}}.
\newblock Idempotent interval analysis and optimization problems.
\newblock {\em Reliable Computing}, 7(5):353--377, 2001.
\newblock E-print \arxiv{math.SC/0101080}.


\bibitem{Lor:93}
M.~Lorenz.
\newblock {\em Object Oriented Software: a Practical Guide}.
\newblock Prentice Hall Books, Englewood Cliffs, N.J., 1993.

\bibitem{MT:03}
G.~G. Magaril-Il'yaev and V.~M. Tikhomirov.
\newblock {\em Convex Analysis: Theory and Applications}, volume 222 of {\em
  Translations of Mathematical Monographs}.
\newblock AMS, Providence, 2003.

\bibitem{Mas:87}
V.~P. Maslov.
\newblock {\em M{\'{e}}thods Op{\'{e}}ratorielles}.
\newblock {\'{E}}ditions MIR, Moscow, 1987.

\bibitem{Mas-86}
V.~P. Maslov.
\newblock New superposition principle for optimization calculus.
\newblock In {\em {S}eminaire sur les {E}quations aux {D}{\'e}riv{\'e}es
  {P}artielles 1985/1986}, {C}entre {M}ath. de l'{\'{E}}cole {P}olytechnique,
  Palaiseau, 1986.
\newblock expos{\'{e}} 24.

\bibitem{Mas-87a}
V.~P. Maslov.
\newblock A new approach to generalized solutions of nonlinear systems.
\newblock {\em Soviet Math. Dokl.}, 42(1):29--33, 1987.

\bibitem{Mas-87b}
V.`P. Maslov.
\newblock On a new superposition principle for optimization problems.
\newblock {\em Uspekhi Math. Nauk [Russian Math. Surveys]}, 42(3):39--48, 1987.

\bibitem{Mas-tech91}
V.~P. Maslov et~al.
\newblock Mathematics of Semirings and its Applications.
\newblock Technical report (in Russian). {I}nstitute for {N}ew {T}echnologies,
  {M}oscow, 1991.

\bibitem{Mas-07}
V.~P. Maslov.
\newblock A general notion of topological spaces of negative dimension and
  quantization of their densities.
\newblock {\em Math. Notes}, 81(1):157--160, 2007.

\bibitem{MS:92}
V.~P. Maslov and S.~N. Samborski{\u{\i}}, editors.
\newblock {\em {Idempotent Analysis}}, volume~13 of {\em Advances in Soviet
  Math.}, Amer. Math. Soc., Providence, 1992.

\bibitem{MV:88}
V.~P. Maslov and K.~A. Volosov, editors.
\newblock {\em Mathematical Aspects of Computer Media}.
\newblock Mir publishers, Moscow, 1988.

\bibitem{Matij}
Yu.~V. Matijasevich.
\newblock A posteriori version of interval analysis.
\newblock In {\em Topics in the {T}heoretical {B}asis and {A}pplications of
  {C}omputer {S}ciences. Proc. of the 4th {H}ung. {C}omp. {S}ci. {C}onf.},
  pages 339--349, Budapest, Akad. Kiado, 1986.

\bibitem{McE:10}
W.~M. McEneaney.
\newblock {\em Max-plus Methods for Nonlinear Control and Estimation}.
\newblock Birkh{\"a}user, Boston et al., 2010.

\bibitem{Mik-05}
G.~Mikhalkin.
\newblock Enumerative tropical algebraic geometry in $\R^2$.
\newblock {\em J. of the ACM}, 18:313--377, 2005.
\newblock E-print \arxiv{math.AG/0312530}.

\bibitem{Mik-06}
G.~Mikhalkin.
\newblock Tropical geometry and its applications.
\newblock In {\em Proceedings of the ICM, volume 2}, Madrid, Spain, pages 827--852, 2006.
\newblock E-print \arxiv{math.AG/0601041}.

\bibitem{Moo:79}
R.~E. Moore.
\newblock {\em Methods and Applications of Interval Analysis}.
\newblock SIAM Studies in Applied Mathematics. SIAM, Philadelphia, 1979.

\bibitem{Mys-05}
H.~My{\v{s}}kova.
\newblock Interval systems of max-separable linear equations.
\newblock {\em Linear Alg. Appl.}, 403:263--272, 2005.

\bibitem{Mys-06}
H.~My{\v{s}}kova.
\newblock Control solvability of interval systems of max-separable linear
  equations.
\newblock {\em Linear Alg. Appl.}, 416:215--223, 2006.

\bibitem{Neu:90}
A.~Neumaier.
\newblock {\em Interval methods for systems of equations}.
\newblock Cambridge University Press, Cambridge, 1990.

\bibitem{Owe-08}
J.~D. Owens.
\newblock GPU computing.
\newblock {\em Proc. of the IEEE}, 96(5):879--899, 2008.

\bibitem{Pan-61}
S.N.N. Pandit.
\newblock A new matrix calculus.
\newblock {\em SIAM J. Appl. Math.}, 9:632--639, 1961.

\bibitem{Pohl:97}
I.~Pohl.
\newblock {\em {O}bject-{O}riented {P}rogramming {U}sing $C++$}.
\newblock Addison-Wesley, Reading, 1997.
\newblock 2nd ed.

\bibitem{Qua-90}
J.~P. Quadrat.
\newblock Th{\'e}or{\`e}ms asymptotiques en programmation dynamique.
\newblock {\em Comptes Rendus Acad. Sci. Paris}, 311:745--748, 1990.

\bibitem{MPlus-Scilab}
J.~P. Quadrat and Max plus~working group.
\newblock Max-plus algebra software.
\newblock \url{http://maxplus.org}; \url{http://scilab.org/contrib};
  \url{http://amadeus.inria.fr},  2007.

\bibitem{RT-87}
Y.~Robert and D.~Tristram.
\newblock An orthogonal systolic array for the algebraic path problem.
\newblock {\em Computing}, 39:187--199, 1987.

\bibitem{Rot-85}
G.~Rote.
\newblock A systolic array algorithm for the algebraic path problem.
\newblock {\em Computing}, 34:191--219, 1985.

\bibitem{Roub}
I.~V. Roublev.
\newblock On minimax and idempotent generalized weak solutions to the
  {H}amilton-{J}acobi equation.
\newblock In G.~L. Litvinov and V.~P. Maslov, editors, {\em Idempotent
  mathematics and mathematical physics}, volume 377 of {\em Contemporary
  Mathematics}, pages 319--338. Amer., Math. Soc., Providence, 2005.

\bibitem{Scha}
H.~H. Schaefer.
\newblock {\em Topological Vector Spaces}.
\newblock Macmillan, New York et al., 1966.

\bibitem{Sed-92}
S.~G. Sedukhin.
\newblock Design and analysis of systolic algorithms for the algebraic path
  problem.
\newblock {\em Computers and artificial intelligence}, 11(3):269--292, 1992.

\bibitem{Shub-92}
M.~A. Shubin.
\newblock Algebraic remarks on idempotent semirings and the kernel theorem in
  spaces of bounded functions.
\newblock In V.~P. Maslov and S.N. Samborski{\u{\i}}, editors, {\em Idempotent
  Analysis}, volume~13 of {\em Advances in Soviet Math.}, pages 151--166, Amer. Math. Soc.,
  Providence, 1992.

\bibitem{Sim-88}
I.~Simon.
\newblock Recognizable sets with multiplicities in the tropical semiring.
\newblock In {\em Lecture Notes in Comp. Sci.}, volume 324, pages 107--120.
Springer, 1988.

\bibitem{Sob-99}
A.~N. Sobolevski{\u{\i}}.
\newblock Interval arithmetic and linear algebra over idempotent semirings.
\newblock {\em Doklady Akademii Nauk}, 369(6):747--749, 1999 (in Russian).
Engl. version: {\em Doklady Math.}, 60(3):431--433, 1999.

\bibitem{SL:94}
A.~Stepanov and M.~Lee.
\newblock {\em The Standard Template Library}.
\newblock Hewlett-Packard, Palo Alto, CA, 1994.

\bibitem{Sub:95}
A.~I. Subbotin.
\newblock {\em Generalized Solutions of First Order {PDE}'s: The Dynamical
  Optimization Perspectives}.
\newblock Birkh{\"a}user, Boston et al., 1995.

\bibitem{Sub-96}
A.~I. Subbotin.
\newblock Minimax solutions of first order partial differential equations.
\newblock {\em Russian Math. Surveys}, 51(2):283--313, 1996.

\bibitem{Vir-00}
O.~Viro.
\newblock {Dequantization of real algebraic geometry on logarithmic paper}.
\newblock In {\em 3rd European Congress of Mathematics: Barcelona, July 10-14,
  2000}, page 135. Birkh{\"a}user, 2001.
\newblock E-print \arxiv{math/0005163}.

\bibitem{Vir-02}
O.~Viro.
\newblock What is an amoeba?
\newblock {\em Notices of the Amer. Math. Soc.}, 49:916--917, 2002.

\bibitem{Vir-08}
O.~Viro.
\newblock From the sixteenth Hilbert problem to tropical geometry.
\newblock {\em Japan. J. Math.}, 3:1--30, 2008.

\bibitem{Voevod}
V.~V. Voevodin.
\newblock {\em Mathematical Foundations of Parallel Computings}.
\newblock {W}orld {S}cientific {P}ubl. {C}o., Singapore, 1992.

\bibitem{Vor-63}
N.~N. Vorobjev.
\newblock The extremal matrix algebra.
\newblock {\em Soviet Math. Dokl.}, 4:1220--1223, 1963.

\bibitem{Vor-67}
N.~N. Vorobjev.
\newblock Extremal algebra of positive matrices.
\newblock {\em Elektron. Informationsverarb. und Kybernetik}, 3:39--71, 1967.

\bibitem{Vor-70}
N.~N. Vorobjev.
\newblock Extremal algebra of nonnegative matrices.
\newblock {\em Elektron. Informationsverarb. und Kybernetik}, 3:302--312, 1970.

\bibitem{Zim-06}
K.~Zimmermann.
\newblock {\em Interval linear systems and optimization problems over
  max-algebras}.  In M.~Fiedler, J.~Nedoma, J.~Ram{\'{\i}}k, J.~Rohn, and K.~Zimmermann.
\newblock {\em Linear Optimization Problems with Inexact Data}.
\newblock Springer, New York, 2006, chapter~6.

\bibitem{IEEE-09}
2009 IEEE International Symposium on Parallel \& Distributed Processing. Rome, Italy,
May 23 -- May 29. ISBN: 978-1-4244-3751-1

\bibitem{ATLAS}
ATLAS: \url{http://math-atlas.sourceforge.net/}

\bibitem{LAPACK}
LAPACK: \url{http://www.netlib.org/lapack/}

\bibitem{PLASMA}
PLASMA: \url{http://icl.cs.utk.edu/plasma/}

\end{thebibliography}

\end{document}